\documentclass[a4paper, 11pt, oneside]{article}

\usepackage{graphicx}
\usepackage[english]{babel}
\usepackage{amsmath}
\usepackage{amssymb}
\usepackage{color}
\usepackage[T1]{fontenc}
\usepackage[utf8]{inputenc}
\usepackage{enumerate}
\usepackage[round]{natbib}
\usepackage{setspace}
\usepackage{floatflt}

\usepackage{nomencl}
\usepackage{ifthen}

\usepackage{geometry}
 \newenvironment{example}{\it{\textbf{Example}}}

\newtheorem{satz}{\bf{Satz}}[section]

\newtheorem{lemma}[satz]{\bf{Lemma}}
\newtheorem{definition}[satz]{\bf{Definition}}

\newtheorem{theorem}[satz]{\bf{Theorem}}
\newtheorem{remark}[satz]{\bf{Remark}}

\newtheorem{corollary}[satz]{\bf{Corollary}}

\newcommand{\R}{\mathbb{R}}
\newcommand{\N}{\mathbb{N}}
\newcommand{\E}{\mathbb{E}}
\newcommand{\Pn}{\mathbb{P}}

\newcommand{\Qn}{\mathbb{Q}}

\newcommand{\F}{\mathcal{F}}

\newcommand{\sZ}{\mathcal{Z}}
\newcommand{\B}{\mathcal{B}}
\newcommand{\A}{\mathcal{A}}

\newcommand{\M}{\mathcal{M}}
\newcommand{\sL}{\mathcal{L}}
\newcommand{\X}{\mathcal{X}}
\newcommand{\Y}{\mathcal{Y}}

\newcommand{\D}{\mathcal{D}}
\newcommand{\W}{\mathbf{W}}
\newcommand{\w}{\mathbf{w}}
\newcommand{\SO}{\mathcal{O}}
\newcommand{\ttt}{\mathbf{t}}

\newcommand{\PN}{P_{\N}}
\newcommand{\QN}{Q_{\N}}
\newcommand{\ZN}{\sZ^{\N}}
\newcommand{\KN}{K_{\N}}
\newcommand{\KNT}{\tilde{K}_{ \N}}

\newcommand{\Summen}{\sum_{i=1}^{n}}
\newcommand{\Z}{$(Z_{i})_{i\in\N}$ }

\newcommand{\Sch}{$(S_{n})_{n\in\N}$ }

\newcommand{\dBL}{d_{\textup{BL}}}

\newcommand{\C}{\mathit{C}}

\begin{document}
\title{Qualitative robustness for bootstrap approximations}
\author{Katharina Strohriegl,University of Bayreuth\\
katharina.strohriegl@uni-bayreuth.de}
\date{\today}
\maketitle

\begin{abstract}
An important property of statistical estimators is qualitative robustness, that is small changes in the distribution of the data only result in small chances of the distribution of the estimator. Moreover, in practice, the distribution of the data is commonly unknown, therefore bootstrap approximations can be used to approximate the distribution of the estimator. Hence qualitative robustness of the statistical estimator under the bootstrap approximation is a desirable property. Currently most theoretical investigations on qualitative robustness assume independent and identically distributed pairs of random variables. However, in practice this assumption is not fulfilled. Therefore, we examine the qualitative robustness of bootstrap approximations for non-i.i.d.\,random variables, for example $\alpha$-mixing and weakly dependent processes. In the i.i.d.\,case qualitative robustness is ensured via the continuity of the statistical operator, representing the estimator, see \cite{Hampel1971} and \cite{CuevasRomo1993}. We show, that qualitative robustness of the bootstrap approximation is still ensured under the assumption that the statistical operator is continuous and under an additional assumption on the stochastic process. In particular, we require a convergence condition of the empirical measure of the underlying process, the so called Varadarajan property.
\end{abstract}

\textbf{AMS:} 60G20, 62G08, 62G09, 62G35

\section{Introduction}

The overwhelming part of theoretical publications in statistical machine learning was done under the assumption that the data is generated by independent and identically distributed (i.i.d.) random variables. However, this assumption is not fulfilled in many practical applications so that non-i.i.d.\ cases increasingly attract attention in machine learning. 
An important property of an estimator is robustness. It is well known that many classical estimators are not robust, which means that small changes in the distribution of the data generating process may highly affect the results, see for example \cite{Huber1981}, \cite{Hampel1968}, \cite{JureckovaPicek2006} or \cite{MaronnaMartin2006} for some books on robust statistics. 
Qualitative robustness is a continuity property of the estimator and means roughly speaking: small changes in the distribution of the data only lead to small changes in the distribution (i.e. the performance) of the estimator. In this way the following kinds of "small errors" are covered: small errors in all data points (rounding errors) and large errors in only a small fraction of the data points (gross errors, outliers). 
Qualitative robustness of estimators has been defined originally in 
\cite{Hampel1968} and \cite{Hampel1971} in the i.i.d.\ case and has been
generalized to estimators for stochastic processes in various ways, for example,
in \cite{PapantoniGray1979}, \cite{Bustos1980}, which will be the one used here, \cite{Cox1981}, \cite{BoenteFairmanYohai1987}, \cite{Zaehle2012}, and \cite{Zaehle2016}, for a more local consideration of qualitative robustness, see for example \cite{Zaehle2017}.

Often the finite sample distribution of the estimator or of the stochastic process of interest is unknown, hence an approximation of the distribution is needed. Commonly, the bootstrap is used to receive an approximation of the unknown finite sample distribution by resampling from the given sample. 

 The classical bootstrap, also called the empirical bootstrap, has been introduced by \cite{Efron1979} for i.i.d.\ random variables.  This concept is based on drawing a bootstrap sample $(Z_1^*,\ldots,Z_m^*)$ of size $m\in\N$ with replacement out of the original sample $(Z_1,\ldots,Z_n)$, $n\in\N$, and approximate the theoretical distribution  $P_n$ of $(Z_1,\ldots,Z_n)$ using the bootstrap sample. For the empirical bootstrap the approximation of the distribution via the bootstrap is given by the empirical distribution of the bootstrap sample $(Z^*_1,\ldots, Z^*_m)$, hence $P^*_n=\otimes_{i=1}^{n}\left(\frac{1}{m} \sum_{i=1}^{m}\delta_{Z_i^*}\right)$, where $\delta_{Z_i}$ denotes the dirac measure. The bootstrap sample itself has distribution $\otimes_{i=1}^{m}\left(\frac{1}{n} \sum_{i=1}^{n}\delta_{Z_i}\right)$.
 
 For an introduction to the bootstrap see for example \cite{Efron1993} and \citet[Chapter 3.6]{VanderVaart1988}. Besides the empirical bootstrap many other bootstrap methods have been developed in order to find good approximations also for non-i.i.d.\,observations, see for example \cite{Kuensch1981}, \cite{Lahiri2003}, and the references therein. In Section \ref{subsec:bootstrapbetamixing} the moving block bootstrap introduced by \cite{Kuensch1989} and \cite{LiuSingh1992} is used to approximate the distribution of an $\alpha$-mixing stochastic process.

It is, also in the non-i.i.d.\,case, still desirable that the estimator is qualitatively robust even for the bootstrap approximation.
That is, the distribution of the estimator under the bootstrap approximation $\sL_{P_n^*}(S_n)$, $n\in\N$, of the assumed, ideal distribution $P_n$ should still be close to the distribution of the estimator under the bootstrap approximation $\sL_{Q _n^*}(S_n)$, $n\in\N$, of the real contaminated distribution $Q_n$. Remember that this is a random object as $P_n^*$ respectively $Q^*_n$ are random.  
 For notational convenience all bootstrap values are noted as usual with an asterisk. 
 
 To show qualitative robustness often generalizations of Hampel's theorem are used, as it is often hard to show qualitative robustness directly. For the i.i.d.\,case Hampel's Theorem ensures qualitative robustness of a sequence of estimators, if these estimators are continuous and can be represented by a statistical operator which is continuous in the distribution of the data generating stochastic process.  Accordingly we try to find results similar to Hampel's theorem for the case of bootstrap approximations for non-i.i.d. cases.
 
Generalizations of Hampel's theorem to  non-i.i.d. cases can be found in \cite{Zaehle2012} and \cite{Zaehle2016}. 
For a slightly different generalization of qualitative robustness, Hampel's theorem has been formulated for strongly stationary and ergodic processes in \cite{Cox1981} and \cite{BoenteFairmanYohai1982}. In \cite{StrohrieglHable2016} a generalization of Hampel's Theorem to a broad class of non-i.i.d.\,stochastic processes is given. \cite{CuevasRomo1993} describes a concept of qualitative robustness of bootstrap approximations for the i.i.d.\,case and for real valued estimators. Also a generalization of Hampel's theorem to this case is given. In \cite{ChristmannVanAelst2013,Christmann2011} qualitative robustness of Efron's bootstrap approximation is shown for the i.i.d.\,case for a class of regularized kernel based learning methods, i.\,e. not necessarily real valued estimators. Moreover \cite{BeutnerZaehle2016} describes consistency of the bootstrap for plug in estimators.

The next chapter contains a definition of qualitative robustness of the bootstrap approximation of an estimator and the main results. In Chapter \ref{subsec:bootstrapid} Theorem \ref{THbootstrapindependentnotidentically} shows qualitative robustness of the bootstrap approximation of an estimator for independent but not necessarily identically distributed random variables, Chapter \ref{subsec:bootstrapbetamixing} contains Theorem \ref{THbootstrapqualitativerobustnessmixing} and \ref{THbootstraphigherdimensions} which generalize the result in \cite{ChristmannVanAelst2013} to $\alpha$-mixing sequences with values in $\R^d$. All proofs are deferred to the appendix.

\section{Qualitative robustness for bootstrap estimators}\label{chapterqualitativerobustnessbootstrap}

Throughout this paper, let $(\sZ,d_{\sZ})$ be a Polish space with some metric $d_{\sZ}$ and Borel-$\sigma$-algebra $\mathcal{B}$. Denote by $\M(\sZ^{\N})$ the set of all probability measures on $(\mathcal{Z}^{\N},\B^{\otimes \N})$. Let 
$(\sZ^{\N},\B^{\otimes\N},\M(\sZ^{\N}))$ be the underlying statistical model.
If nothing else is stated, we always use Borel-$\sigma$-algebras for topological spaces. 
Let $(Z_i)_{i\in\N}$ be the coordinate process on $\sZ^{\N}$, that is $Z_i:\ZN\rightarrow\sZ,\;(z_j)_{j\in\N}\mapsto z_i$, $i\in\N$. Then the process has law $\PN$ under $\PN\in\M(\ZN)$. Moreover let $P_n:=(Z_1,\ldots, Z_n)\circ\PN$ be the $n$-th order marginal distribution of $\PN$ for every $n\in\N$ and $\PN\in\M(\ZN)$.
We are concerned with a sequence of estimators
$(S_n)_{n\in\mathbb{N}}$ on the stochastic process $(Z_i)_{i\in\mathbb{N}}$.
The estimator may take its values in any Polish space $H$ with some metric $d_H$; that is,
$S_n:\mathcal{Z}^n\rightarrow H$ for every $n\in\mathbb{N}$.

Our work applies to estimators which can be
represented by a statistical operator 
$S:\mathcal{M}(\mathcal{Z})\rightarrow H$, that is,
\begin{eqnarray}\label{represented-by-statistical-functional}
  S\big(\Pn_{\w_n}\big)
  \;=\;S_n(\w_n)\;=\;S_n(z_1,\dots,z_n)
  \qquad\forall\,\w_n=(z_1,\dots,z_n)\in\mathcal{Z}^n,
  \quad\forall\,n\in\mathbb{N},
\end{eqnarray}

where $\Pn_{\w_n}$ denotes the empirical measure defined
by $\Pn_{\w_n}(B):=\frac{1}{n}\Summen I_B(z_i)$, $B\in\mathcal{B}$,
for the 
observations $\w_n=(z_1,...,z_n)\in\mathcal{Z}^n$.
Examples of such estimators are M-estimators, R-estimators, see \citet[Theorem 2.6]{Huber1981}, or Support Vector Machines, see \cite{HableChristmann2011}.


Based on the generalization of Hampel's concept of $\Pi$-robustness from \cite{Bustos1980}, we define qualitative robustness for bootstrap approximations for non-i.i.d sequences of random variables. The stronger concept of $\Pi$-robustness is needed here, as we do not assume to have i.i.d.\,random variables, which are used in \cite{CuevasRomo1993}.

Therefore the definition of qualitative robustness stated below is stronger than the definition in \cite{CuevasRomo1993}, i.\,e. if we use this definition for the i.i.d.\,case the assumption $\dBL(P_n,Q_n)=\dBL(\otimes_{i=1}^n P, \otimes_{i=1}^n Q)<\delta$ implies $\dBL(P,Q)<\delta$, where $\dBL$ denotes the bounded Lipschitz metric. This can be seen similar to the proof of Lemma \ref{Lemma: Zusammenhang Mischmass und Produktmass} in Section \ref{subsec:bootstrapid}.

Now, let $P_{\N}^*$ be the approximation of $P_{\N}$ with respect to the bootstrap. Define the bootstrap sample $(Z^*_1,\ldots,Z^*_n)$ as the first $n$ coordinate projections $Z^*_i:\ZN\rightarrow\sZ$, where the law of the stochastic process $(Z^*_i)_{i\in\N}$ has to be chosen according to the bootstrap procedure. For the empirical bootstrap, for example, the bootstrap sample is chosen via drawing with replacement from the given observations $z_1,\ldots , z_{\ell}$, $\ell\in\N$. Hence the distribution of the bootstrap sample is $\otimes_{n\in\N} \frac{1}{\ell}\sum_{i=1}^{\ell} \delta_{z_i}$, with finite sample distributions $\otimes_{j=1}^n\frac{1}{\ell}\sum_{i=1}^{\ell}\delta_{z_i}=(Z^*_1,\ldots,Z^*_n)\left(\otimes_{n\in\N} \frac{1}{\ell}\sum_{i=1}^{\ell} \delta_{z_i}\right)$.\\
Contrarily to the classical case of qualitative robustness the distribution of the estimator under $P^*_n$, $\sL_{P^*_n}(S_n)$ is a random probability measure, as the distribution $P^*_n=\otimes_{i=1}^n \frac{1}{\ell}\sum_{i=1}^{\ell} \delta_{Z^*_i}$, $Z^*_i:\ZN\rightarrow \sZ$, is random. Hence the mapping $z_{\N}\mapsto\sL_{P^*_n}(S_n)$, $z_{\N}\in\ZN$, is itself a random variable with values in $\M(H)$, i.\,e. on the space of probability measures on $H$, equipped with the weak topology on $\M(H)$. The measurability of this mapping is ensured by \citet[Lemma D1]{BeutnerZaehle2016}. 

Contrarily to the original definitions of qualitative robustness in \cite{Bustos1980} the bounded Lipschitz metric $\dBL$ is used  instead of the Prohorov metric $\pi$ for the definition of qualitative robustness of the bootstrap approximation below. This is equivalent to \cite{CuevasRomo1993}. Let $\X$ be a separable metric space, then the bounded Lipschitz metric on the space of probability measures $\M(\X)$ on $\X$ is defined by:
	\[
	d_{\textup{BL}}(P,Q):=\sup\left\{\left|\int f \,dP-\int f\,dQ\right|;\; f\in \textup{BL}(\X), \|f\|_{\textup{BL}}\leq 1\right\}
	\]
where $\|\cdot\|_{\textup{BL}}:=|\cdot|_1+\|\cdot\|_{\infty}$ denotes the bounded Lipschitz norm with $|f|_{1}=\sup_{x\neq y}\frac{|f(x)-f(y)|}{d(x,y)}$ and $\|\cdot\|_{\infty}$ the supremum norm $\|f\|_{\infty}:=\sup_{x}|f(x)|$ and the space of bounded Lipschitz functions is defined as $\textup{BL}:=\{f:\X\rightarrow\R\;|\; f\;\text{Lipshitz and } \|f\|_{\textup{BL}}<\infty\}$. This is due to technical reasons only. Both metrics metricize the weak topology on the space of all probability measures $\M(\X)$, for Polish spaces $\X$, see, for example, \citet[Chapter 2, Corollary 4.3]{Huber1981} or \citet[Theorem 11.3.3]{Dudley1989}, and therefore can be replaced while adapting $\delta$ on the left hand-side of implication \eqref{robustnessbootstrap}. 
If $\X$ is a Polish space, so is $\M(\X)$ with respect to the weak topology, see \citet[Chapter 2, Theorem 3.9]{Huber1981}. 
Hence the bounded Lipschitz metric on the right-hand side of implication \eqref{robustnessbootstrap} operates on a space of probability measures on the Polish space $\M(\X)$. Therefore the Prohorov metric and the bounded Lipschitz metric can again be replaced while adapting $\varepsilon$ in \eqref{robustnessbootstrap}. Similar to \cite{CuevasRomo1993} the proof of the theorems below rely on the fact that the set of bounded Lipschitz functions $\textup{BL}$ is a uniform Glivenko-Cantelli class, which implies uniform convergence of the bounded Lipschitz metric of the empirical measure to a limiting distribution, see \cite{DudleyZinn1989}. Therefore the definition is given with respect to the bounded Lipschitz metric.

		\begin{definition}[Qualitative robustness for bootstrap approximations]\label{Defqualitativerobustnessbootstrap}\hspace{0,5cm}\\
		Let $P_{\N}\in\M(\sZ^{\N})$ and let $P^*_{\N}\in\M(\sZ^{\N})$ be the bootstrap approximation of $P_{\N}$. Let $\mathcal{P}\subset\M(\sZ^{\N})$ with $P_{\N}\in\mathcal{P}$. Let
		$S_n:\sZ^{n}\rightarrow H$, $n\in\N$, be a sequence of estimators. 
		Then the sequence of bootstrap approximations $(\sL_{P_n^*}(S_n))_{n\in\N}$ is called qualitatively robust at $\PN$ with respect to $\mathcal{P}$
		if, for every $\varepsilon>0$, there is $\delta>0$ 
		such that there is $n_0\in\N$ such that for every $n\geq n_0$ and for every $\QN\in\mathcal{P}$,	
		\begin{equation}\label{robustnessbootstrap}
		d_{\textup{BL}}(P_n,Q_n)<\delta\;
		\Rightarrow\;
		d_{\textup{BL}}(\sL(\sL_{P_n^*}(S_n)),\sL(\sL_{Q_n^*}(S_n))) < \varepsilon.
		\end{equation}
		
	Here $\sL(\sL_{P_n^*}(S_n))$ (respectively $\sL(\sL_{Q_n^*}(S_n))$) denotes the distribution of the bootstrap approximation 
	of the estimator $S_n$ under $P^*_n$ (respectively $Q^*_n$).
		\end{definition}

This definition of qualitative robustness with respect to the subset $\mathcal{P}$ indicates that we do not show \eqref{robustnessbootstrap} for arbitrary probability measures $\QN\in\M(\sZ^{\N})$. All of our results require the contaminated process to at least have the same structure as the ideal process. 
This is due to the use of the bootstrap procedure. The empirical bootstrap, which is used below, only works well for a few processes, see for example \cite{Lahiri2003}, hence the assumptions on the contaminated process are necessary. To our best knowledge there are no results concerning qualitative robustness of the bootstrap approximation for general stochastic processes without any assumptions on the second process and it is probably very hard to show this for every $\QN\in\M(\sZ^{\N})$, respectively $\mathcal{P}=\M(\sZ^{\N})$.
Another difference to the classical definition of qualitative robustness in \cite{Bustos1980} is the restriction to $n\geq n_0$. As the results for the bootstrap are asymptotic results, we can not achieve the equicontinuity for every $n\in\N$, but only asymptotically.

As the estimators can be represented by a statistical operator which depends on the empirical measure it is crucial to concern stochastic processes which at last provide convergence of their empirical measure. Therefore, \cite{StrohrieglHable2016} proposed to choose Varadarajan process. Let $(\Omega,\A,\mu)$ be a probability space. Let $(Z_i)_{i\in\N},$ $Z_i:\Omega\rightarrow\sZ$, $i\in\N$, be a stochastic process and $\W_n:=(Z_1,\ldots, Z_n)$. Then the stochastic process $(Z_i)_{i\in\N}$ is called a \emph{(strong) 
  Varadarajan process} if there exists a probability measure $P\in\M(\sZ)$ such that 
  \begin{equation*}
	\pi(\Pn_{\W_n},P)\;\xrightarrow[\;n\rightarrow\infty\;]{}\;0
    \quad\text{almost surely.}
  \end{equation*}
  The stochastic process \Z is called \emph{weak Varadarajan process} if
  \begin{equation*}
	 \pi(\Pn_{\W_n},P)\;\xrightarrow[\;n\rightarrow\infty\;]{}\;0
     \quad\text{in probability.}
  \end{equation*}
 
Examples for Varadarajan processes are certain Markov Chains, some mixing processes, ergodic process and processes which satisfy a law of large numbers for events in the sense of \citet[Definition 2.1]{SteinwartHush2009}, see \cite{StrohrieglHable2016} for details.

 \subsection{Qualitative robustness for independent not identically distributed processes}\label{subsec:bootstrapid}

In this section we relax the i.i.d.\,assumption in view of the identical distribution. We assume the random variables $Z_i$, $i\in\N$, to be independent, but not necessarily identically distributed. \\
The result below generalizes \citet[Theorem 3]{ChristmannVanAelst2013} and \cite{Christmann2011}, as the assumptions on the stochastic process are weaker as well as those on the statistical operator. Compared to Theorem 3 in \citet{CuevasRomo1993}, which shows qualitative robustness of the sequence of bootstrap estimators with values in $\R$, we have to strengthen the assumptions on the sample space, but do not need the estimator to be uniformly continuous. But keep in mind, that the assumption $\dBL(P_n,Q_n)< \delta$ implies $\dBL(P,Q)<\delta$, which is used for the i.i.d. case, in \cite{ChristmannVanAelst2013} and \cite{CuevasRomo1993}.
%
%
%
%
%
%

		\begin{theorem}\label{THbootstrapindependentnotidentically}
		Let the sequence of estimators \Sch be represented by a statistical operator 
  		$S:\M(\sZ)\rightarrow H$
  		via (\ref{represented-by-statistical-functional}) for a Polish space $H$ and let $(\sZ,d_{\sZ})$ be a totally bounded metric space. \\		
  		Let $\PN=\otimes_{i\in\N}P^i$, $P^i\in\M(\sZ)$ be an infinite product measure such that the coordinate process $(Z_i)_{i\in\N}$, $Z_i\colon \ZN\rightarrow z_i$, $i\in\N$, is a strong Varadarajan process with limiting distribution $P$. Moreover define $\mathcal{P}:=\left\{\QN\in\M(\sZ^{\N});\; \QN=\otimes_{i\in\N}Q^i,\;Q^i\in\M(\sZ)\right\}$. Let $S:\mathcal{M}(\mathcal{Z}) \rightarrow H$ be continuous at $P$ with respect to $\dBL$
  		and let the estimators $S_n: \sZ^{n}\rightarrow H,\;n\in\N$, be continuous.\\  		
  		Then the sequence of bootstrap approximations $(\sL_{P_n^*}(S_n))_{n\in\N}$, is qualitatively robust at $\PN$ with respect to $\mathcal{P}$.
		\end{theorem}

\begin{remark}
	The required properties on the statistical operator $S$ and on the sequence of estimators \Sch 
	in Theorem \ref{THbootstrapindependentnotidentically} ensure the qualitative robustness of $(S_n)_{n\in\N}$, 
	as long as the assumptions on the underlying stochastic processes are fulfilled. \\
	The proof shows that the bootstrap approximation 
	of every sequence of estimators \Sch which is qualitatively robust in the sense of the definitions in \cite{Bustos1980} and \citet[Definition 1]{StrohrieglHable2016} is qualitatively robust in the sense of Theorem \ref{THbootstrapindependentnotidentically}.
	\end{remark}

Hence Hampel's theorem for the i.i.d.\,case can be generalized to bootstrap approximations and to the case of not necessarily identically distributed random variables if qualitative robustness is based on the definition of $\Pi$-robustness.

 Unfortunately, the assumption on the space $(\sZ,d_{\sZ})$ to be totally bounded seems to be necessary. In the proof of Theorem \ref{THbootstrapindependentnotidentically} we use a result of \cite{DudleyZinn1989} to show uniformity on the space of probability measures $\M(\sZ)$. This result needs the bounded Lipschitz functions to be a uniform Glivenko-Cantelli class, which is equivalent to $(\sZ,d_{\sZ})$ being totally bounded, see \citet[Proposition 12]{DudleyZinn1989}. In order to weaken the assumption on $(\sZ,d_{\sZ})$, probably another way to show uniformity on the space of probability measures $\M(\sZ)$ has to be found.





	A short look on the metrics  used on $\sZ^n$ is advisable. We consider $\mathcal{Z}^n$ as the 
$n$-fold product space of the Polish space $(\mathcal{Z},d_{\sZ})$. 
The product space $\mathcal{Z}^n$ is again a Polish space (in the product
topology) and it is tempting to use a $p$-product metric $d_{n,p}$ on
$\mathcal{Z}^n$, that is,
	\begin{equation}
			d_{n,p}\big((z_1,\dots,z_n),(z_1^\prime,\dots,z_n^\prime)\big)
 			 \;=\;\big\|\big(d_{\sZ}(z_1,z_1^\prime),\dots,d_{\sZ}(z_n,z_n^\prime)\big)\big\|_{p}\label{pnNorm}
	\end{equation}
where $\|\cdot\|_p$ is a $p_n$-norm on $\mathbb{R}^n$ for
$1\leq p\leq \infty$.
For example, $d_{n,2}$ is the Euclidean metric on $\R^n$ and
$d_{n,\infty}\big((z_1,\dots,z_n),(z_1^\prime,\dots,z_n^\prime)\big)
  =\max_i d(z_i,z_i^\prime)
$; all these metrics are strongly equivalent.
However, these common metrics do not cover the intuitive meaning of qualitative robustness as the distance between two
points in $\mathcal{Z}^n$ (i.e., two data sets) 
is small only if all coordinates are close together (small rounding errors). So points where only a small fraction of the coordinates are far-off (gross errors) are excluded. Using these metrics, the qualitative robustness of the sample mean at every $\PN\in\M(\sZ^{\N})$ can be shown, see e.g. \citet[Proposition 1]{StrohrieglHable2016}. But the sample mean is a highly non-robust estimator, as gross errors have great impact on the estimate. Following
\cite{BoenteFairmanYohai1987}, we use the metric $d_n$ on $\mathcal{Z}^n:$ \begin{equation}\label{infmetrik}
				d_{n}\big((z_1,\dots,z_n),(z_1^\prime,\dots,z_n^\prime)\big)
 				 \;=\;\inf\big\{\varepsilon>0: 
                 \sharp\{i: d(z_i,z_i^\prime)\geq\varepsilon\}/n\leq\varepsilon
          			 \big\}\,.
		\end{equation}
This metric on $\sZ^{n}$ covers both kinds of "small errors". Though $d_n$ is not strongly equivalent to $d_{n,p}$ in general,
it is topologically equivalent to the $p$-product metrics $d_{n,p}$, see \citet[Lemma 1]{StrohrieglHable2016}. Hence, $\sZ^n$ is metrizable also with metric $d_n$. Moreover the continuity of $S_n$ on $\sZ^n$ is with respect to the product topology on $\sZ^n$ which can, due to the topological equivalence of these two metrics, be seen with respect to the common metrics $d_{n,p}$.



The next part gives two examples of stochastic processes of independent, but not necessarily identically distributed random variables, which are Varadarajan processes. In particular these stochastic processes even satisfy a strong law of large numbers for events (SLLNE) in the sense of \cite{SteinwartHush2009} and therefore are, due to \citet[Theorem 2]{StrohrieglHable2016}, strong Varadarajan processes. The first example is rather simple and describes a sequence of univariate normal distributions.
	
	\begin{example} \textbf{1} 
	Let $(a_i)_{i\in\N}\subset \R$ be a sequence with $\lim_{i\rightarrow \infty} a_i= a\in\R$ and let $|a_i|\leq c$, for some constant $c>0$ for all $i\in \N$. 
	Let $(Z_i)_{i\in\N},$ $Z_i:\Omega\rightarrow \R$, be a stochastic process where $Z_i$, $i\in\N$, are independent and $Z_i\sim N(a_i,1),\; i\in \N$.
	Then the process \Z is a strong Varadarajan process.
	\end{example}

	The second example are stochastic processes where the distributions of the random variables $Z_i$, $i\in\N$,
	are lying in a so-called shrinking $\varepsilon$-neighbourhood of a probability measure $P$.

\begin{example} \textbf{2}
	Let $(\sZ,\B)$ be a measurable space and let \Z be a stochastic process with independent random variables $Z_i:\Omega\rightarrow \sZ$, $Z_i\sim P^i$, where 
	\[
	P^i=(1-\varepsilon_i)P+\varepsilon \tilde{P}^i
	\]
	for a sequence $\varepsilon_i\rightarrow0$, $i\rightarrow\infty$, $\varepsilon_i>0$ and $\tilde{P}^i,\;P\in\M(\sZ)$, $i\in\N$. 
	Then the process \Z is a strong Varadarajan process.
	\end{example}

	The next corollary shows, that Support Vector Machines are qualitatively robust. For a detailed introduction to Support Vector Machines see e.g., \cite{SchoelkopfSmola2002} and \cite{SteinwartChristmann2008}. Let $D_n:=(z_1,z_2,\ldots,z_n)=((x_1,y_1),(x_2,y_2),\ldots,(x_n,y_n))$ be a given dataset.
	
	\begin{corollary}\label{COBootstrapidSVM}
	Let $\sZ=\X\times\Y$, $\Y\subset\R$ closed, be a totally bounded, metric space and let \Z be a stochastic process where the random variables $Z_i$, $i\in\N$, are independent and $Z_i\sim P^i:= (1-\varepsilon_i)P+\varepsilon_i \tilde{P}^i$, $P,\tilde{P}^i\in \M(\sZ)$.	
	Moreover let $(\lambda_n)_{n\in\N}$ be a sequence of positive real valued numbers with $\lambda_n\rightarrow\lambda_0,\;
		  n\rightarrow\infty$, for some $\lambda_0>0$. Let $H$ be a reproducing kernel Hilbert space with continuous and bounded kernel $k$ and let $S_{\lambda_n}: (\X\times\Y)^n\rightarrow H$ be the SVM estimator, 
		  which maps $D_n$ to $f_{L^*,D_n,\lambda_n}$ for a continuous and convex loss function 
		  $L: \X\times\Y\times\Y\rightarrow[0,\infty[$. 
		  It is assumed that $L(x,y,y)=0$ for every $(x,y)\in\X\times\Y$ and that $L$ is additionally 
		  Lipschitz continuous in the last argument.\\		 
	Then we have 
  		for every $\varepsilon >0$ there is $\delta>0$ such that there is $n_0\in\N$ such that for all $n\geq n_0$ 
  		and for every process $(\tilde{Z}_i)_{i\in\N}$,
  		where $\tilde{Z}_i$ are independent and have distribution $Q^i$, $i\in\N$: 
  			\[
  			d_{\textup{BL}}(P_n,Q_n)<\delta\;\Rightarrow\;d_{\textup{BL}}(\sL(\sL_{P_n^*}(S_n)),\sL(\sL_{Q_n^*}(S_n)))<\varepsilon.
  			\]
	\end{corollary}

	That is, the sequence of bootstrap approximations is qualitatively robust if the second (contaminated) process $(\tilde{Z}_i)_{i\in\N}$ is still of the same kind, i.e. still independent, as the original uncontaminated process $(Z_i)_{i\in\N}$.


\subsection{Qualitative robustness for the moving block bootstrap of $\alpha$-mixing processes}\label{subsec:bootstrapbetamixing}

Dropping the independence assumption we now focus on real valued mixing processes, in particular on strongly stationary $\alpha$-mixing or strong mixing stochastic processes. The mixing notion is an often used and well-accepted dependence notion which quantifies the degree of dependence of a stochastic process. There exist several types of mixing coefficients, but all of them are based on differences
between probabilities $\mu(A_1\cap A_2)-\mu(A_1)\mu(A_2)$. There is a large literature on this dependence structure. 
For a detailed overview on mixing, see \cite{Bradley2005}, \cite{Bradley20071,Bradley20072,Bradley20073}, and \cite{Doukhan1994} and the references therein. The $\alpha$-mixing structure has been introduced in \cite{Rosenblatt1956}. Also examples of relations between dependence structures and mixing coefficients can be found in the references above. 
Let $\Omega$ be a set equipped with two $\sigma$-algebras $\A_1$ and $\A_2$ and a probability measure $\mu$. Then the $\alpha$-mixing coefficient is defined by
$$\alpha(\A_1,\A_2,\mu)		 :=\sup\{|\mu(A_1\cap A_2)-\mu(A_2)\mu(A_2)|\; |\; A_1\in\A_1,\; A_2\in\A_2\}.$$
By definition the coefficients equal zero, if the $\sigma$-algebras are independent. 

Moreover mixing can be defined for stochastic processes. We follow \citet[Definition 3.1]{SteinwartHush2009}:
\begin{definition}\label{Def: bootstrap weakly alpha bi mixing}
Let \Z be a stochastic process, $Z_i:\Omega\rightarrow \sZ$, $i\in\N$, and let $\sigma(Z_i)$ be the $\sigma$-algebra generated by $Z_i$, $i\in\N$. Then the $\alpha$-$bi$- and the $\alpha$-mixing coefficients are defined by
\begin{align*}
&\alpha((Z)_{i\in\N},\mu,i,j)=\alpha(\sigma(Z_i),\sigma(Z_j),\mu)\\
&\alpha((Z)_{i\in\N},\mu,n)=\sup_{i\geq 1} \alpha(\sigma(Z_i),\sigma(Z_{i+n}),\mu).
\end{align*}
A stochastic process \Z is called $\alpha$-
mixing with respect to $\mu$ if 
\begin{align*}
\lim_{n\rightarrow \infty} \alpha((Z)_{i\in\N},\mu,n)= 0.
\end{align*}
It is called weakly 
$\alpha$-$bi$-mixing with respect to $\mu$ if 
\begin{align*}
 \lim_{n\rightarrow\infty} \frac{1}{n^2}\sum_{i=1}^n\sum_{j=1}^{n}\alpha((Z)_{i\in\N},\mu,i,j)=0.
\end{align*}
\end{definition}

Instead of Efron's empirical bootstrap another bootstrap approach is used in order to represent the dependence structure of an $\alpha$-mixing process. \cite{Kuensch1989} and \cite{LiuSingh1992} introduced the moving block bootstrap (MBB). Often resampling of single observations can not preserve the dependence structure of the process, therefore they decided to take blocks of length $b$ of observations instead. The dependence structure of the process is preserved, within these blocks. The block length $b$ increases with the number of observations $n$ for asymptotic considerations. A slight modification of the original moving block bootstrap, see for example \cite{Politis1990} and \cite{ShaoYu1993}, is used in the next two theorems in order to avoid edge effects. 

The proofs are based on central limit theorems for empirical processes. There are several results concerning the moving block bootstrap of the empirical process in case of mixing processes, see for example \cite{Buehlmann1994}, \cite{Naik1994}, and \citet[Theorem 2.2]{Peligrad1998} for $\alpha$-mixing sequences and \cite{Radulovic1996} and \cite{Buehlmann1995} for $\beta$-mixing sequences. To our best knowledge there are so far no results concerning qualitative robustness for bootstrap approximations of estimators for $\alpha$-mixing stochastic processes. Therefore, Theorem \ref{THbootstrapqualitativerobustnessmixing} shows qualitative robustness for a stochastic process with values in $\R$. The proof is based on \citet[Theorem 2.2]{Peligrad1998}, which provides a central limit theorem under assumptions on the process, which are weaker than those in \cite{Buehlmann1994} and \cite{Naik1994}. In the case of $\R^d$-valued, $d>1$, stochastic processes, stronger assumptions on the stochastic process are needed, as the central limit theorem in \cite{Buehlmann1994} requires stronger assumptions, see Theorem \ref{THbootstraphigherdimensions}.

Let $Z_1,\ldots, Z_n$, $n\in\N$, be the first $n$ projections of a real valued stochastic process \Z and let $b\in\N, b<n$, be the block length. Then, for fixed $n\in\N$, the sample can be divided into blocks $B_{i,b}:=(Z_i,\ldots, Z_{i+b-1})$. If $i>n-b+1$, we define $Z_{n+j}= Z_j$, for the missing elements of the blocks. To get the MBB bootstrap sample $\W^*_n=(Z^*_1,\ldots, Z^*_n)$, $\ell$ numbers $I_1,\ldots, I_{\ell}$ from the set $\{1,\ldots,n\}$ are randomly chosen with replacement. 
Without loss of generality it is assumed that $n=\ell b$, if $n$ is not a multiple of $b$ we simply cut the last block, which is usually done in literature. 
Then the sample consists of the blocks $B_{I_1,b},B_{I_2,b},\ldots, B_{I_{\ell},b}$, that is $Z^*_1= Z_{I_1},Z^*_2= Z_{I_1+1},\ldots, Z^*_{b}=_{I_1+b-1}, Z^*_{b+1}= Z_{I_2},\ldots,Z^*_{\ell b}=Z_{I_{\ell}+b-1}$. 

As we are interested in estimators $S_n$, $n\in\N$, which can be represented by a statistical operator $S\colon\M(\sZ)\rightarrow H$ via $S(\Pn_{\w_n})=S_n(z_1,\ldots,z_n)$, for a Polish space $H$, see \eqref{represented-by-statistical-functional}, the empirical measure of the bootstrap sample $\Pn_{\W^*_n}=\frac{1}{n}\sum_{i=1}^{n}\delta_{Z^*_i}$ should approximate the empirical measure of the original sample $\Pn_{\W_n}=\frac{1}{n}\sum_{i=1}^{n}\delta_{Z_i}$.
Contrarily to qualitative robustness in the case of independent and not necessarily identically distributed random variables (Theorem \ref{THbootstrapindependentnotidentically}), the assumptions on the statistical operator $S$ are strengthened for the case of $\alpha$-mixing sequences. In particular the statistical operator $S$ is assumed to be uniformly continuous for all $P\in (\M(\sZ),\dBL)$.
For the first theorem we assume the random variables $Z_i,\; i \in \N$, to be real valued and bounded. Without loss of generality we assume $0\leq Z_i\leq 1$, otherwise a transformation leads to this assumption.
 For the bootstrap for the true as well as for the contaminated process, we assume the block length $b(n)$ and the number of blocks $\ell(n)$ to be sequences of integers satisfying
 	$$n^h\in \SO(b(n)),\; b(n)\in \SO(n^{1/3-a}),\;\text{for some}\; 0<h<\frac{1}{3}-a,\;0<a<	\frac{1}{3},$$
 	$b(n)=b(2^{q})$ for $2^{q}\leq n < 2^{q+1}, \; q\in\N,\;b(n)\rightarrow\infty$, $n\rightarrow\infty$ and $b(n)\cdot \ell(n)=n$, $n\in\N$.	

	\begin{theorem}\label{THbootstrapqualitativerobustnessmixing}
	Let $\PN\in\M(\R^{\N})$ be a probability measure on $(\R^{\N},\B^{\otimes\N})$ such that the coordinate process $\left(Z_i\right)_{i\in\N}$, $Z_i:\R^{\N}\rightarrow\R$ is bounded, strongly stationary, and $\alpha$-mixing with 
	\begin{equation}\label{THbootstrap alphamixing bedinung 1d}
	\sum_{m>n}\alpha (\sigma(Z_1,\ldots,Z_i),\sigma(Z_{i+m},\ldots),\PN) = \SO(n^{-\gamma}),\;i\in\N,\;\text{for some}\;\gamma>0.
	\end{equation}	
	Let $\mathcal{P}\subset\M(\R^{\N})$ be the set of probability measures such that the coordinate process fulfils the properties above for the same $\gamma >0$.
 	Let $H$ be a Polish space, with some metric $d_H$, let \Sch  be a sequence of estimators which can be represented by a statistical operator $S:\M(\R)\rightarrow H$ via \eqref{represented-by-statistical-functional}. Moreover let $S_n$ be continuous and let $S$ be additionally uniformly continuous with respect to $\dBL$. 
 	Then the sequence of estimators \Sch is qualitatively robust at $\PN$ with respect to $\mathcal{P}$.
	\end{theorem}

The assumptions on the stochastic process are on the one hand, together with the assumptions on the block length, used to ensure the validity of the bootstrap approximation and on the other hand, together with the assumptions on the statistical operator, respectively the sequence of estimators, to ensure the qualitative robustness. 




The next theorem generalizes this result to stochastic processes with values in $[0,1]^d$, $d>1$, instead of $[0,1]\subset\R$. Therefore, for example, the bootstrap version of the SVM estimator is qualitatively robust under weak conditions.
The proof of the next theorem follows the same lines as the proof of the theorem above, but another central limit theorem, which is shown in \cite{Buehlmann1994}, is used. Therefore the assumptions on the mixing property of the stochastic process are stronger and the random variables $Z_i,\;i\in\N,$ are assumed to have continuous marginal distributions. Again the bootstrap sample results of a moving block bootstrap where $\ell(n)$ blocks of length $b(n)$ are chosen, again assuming $\ell(n)\cdot b(n) =n$. Moreover, let $b(n)$ be a sequences of integers satisfying
 	$$b(n)=\SO(n^{\frac{1}{2}-a}) \;\text{for some}\; a>0.$$

\begin{theorem}\label{THbootstraphigherdimensions}
Assume $\sZ=[0,1]^d$, $d>1$. Let $\PN$ be a probability measure such that the coordinate process $\left(Z_i\right)_{i\in\N}$, $Z_i:\sZ^{\N}\rightarrow \sZ$ is strongly stationary and  $\alpha$-mixing with 
\begin{equation}\label{proof: bootstrapmixing Rd Bedingung an Pozess}
\sum_{m=0}^{\infty}(m+1)^{8d+7}(\alpha (\sigma(Z_1,\ldots,Z_i),\sigma(Z_{i+m},\ldots),\PN))^{\frac{1}{2}}<\infty,\; i\in\N.
\end{equation}
Assume that $Z_i$ has continuous marginal distributions for all $i\in\N$.
Define the set of probability measures $\mathcal{P}\subset \M(\sZ)$ such that the coordinate process is strongly stationary and $\alpha$-mixing as in \eqref{proof: bootstrapmixing Rd Bedingung an Pozess}.\\
 	Let $H$ be a Polish space, wit some metric $d_H$, \Sch be a sequence of estimators such that $S_n:\sZ^n\rightarrow H$ is continuous and assume that 
 	$S_n$ can be represented by a statistical operator $S:\M(\sZ)\rightarrow H$ via \eqref{represented-by-statistical-functional} 
 	which is additionally uniformly continuous with respect to $\dBL$.
 	 
 	Then the sequence of estimators \Sch is qualitatively robust at $\PN$ with respect to $\mathcal{P}$.	
\end{theorem}

Although the assumptions on the statistical operator $S$, compared to Theorem \ref{THbootstrapindependentnotidentically}, were strengthened in order to generalize the qualitative robustness to $\alpha$-mixing sequences in Theorem \ref{THbootstrapqualitativerobustnessmixing} and \ref{THbootstraphigherdimensions}, M-estimators are still an example for qualitative robust estimators if the sample space $(\sZ,d_{\sZ})$, $\sZ\subset \R$ is compact. The compactness of $(\sZ,d_{\sZ})$ implies the compactness of the space $(\M(\sZ),\dBL)$, see \citet[Theorem 6.4]{Parthasarathy1967}. As the statistical operator $S$ is continuous, the compactness of $\M(\sZ)$ implies the uniform continuity of $S$. Another example of M-estimators which are uniformly continuous even if the input space is not compact is given in \citet[Theorem 4]{CuevasRomo1993}. 


	\textbf{Acknowledgements:} This research was partially supported by the DFG Grant 291/2-1 "Support Vector Machines bei stochastischer Unabhängigkeit". Moreover I would like to thank Andreas Christmann for helpful discussions on this topic.

\section{Proofs}
This section contains the proofs of the main theorems and corollaries.

\subsection{ Proofs of Section \ref{subsec:bootstrapid}}

Before proving Theorem \ref{subsec:bootstrapid}, we state a rather technical lemma, connecting the  product measure $\otimes_{i=1}^n P^i\in\M(\sZ^n)$ of independent random variables to their mixture measure $\frac{1}{n}\sum_{i=1}^n P^i\in\M(\sZ)$. Let $(\sZ,d_{\sZ})$ be a Polish space.

	\begin{lemma}\label{Lemma: Zusammenhang Mischmass und Produktmass}
		Let $P_n, Q_n\in\M(\sZ^n)$ such that $P_n=\otimes_{i=1}^n P^i$ and 
		$Q_n=\otimes_{i=1}^n Q^i$, $P^i,Q^i\in\M(\sZ),\; i\in\N$.
		Then for all $\delta >0$:
		\begin{equation*}
			\dBL(P_n,Q_n)\leq \delta\quad\Rightarrow\quad 
			\dBL\left(\frac{1}{n}\sum_{i=1}^n P^i,\frac{1}{n}\sum_{i=1}^n Q^i\right)\leq \delta.
		\end{equation*}
	\end{lemma}


\textbf{Proof:} Let $\textup{BL}_1$ be the set of bounded Lipschitz functions with $\|f\|_{\textup{BL}}\leq 1$.By assumption we have $\dBL(P_n,Q_n)\leq \delta$. Moreover for a function $f:\sZ\rightarrow \R$:
\begin{align}\label{proof:qualitative robustness boottsrap Lemma produktmass und Mischmass}
\int_{\sZ}f(z_i) \,dP^i(z_i)=\int_{\sZ^{n-1}}\int_{\sZ}f(z_i) \,d  P^i(z_i) \,d\left( \otimes_{j\neq i}P^j(z_j)\right).
\end{align} 

Then, 
		\begin{align*}
			 \sup_{f\in \textup{BL}_1(\sZ)} & \left| \int_{\sZ}f(z_i)\; d \left[\frac{1}{n}\sum_{i=1}^n P^i(z_i)\right]-
			 \int_{\sZ}f(z_i)\,d \left[\frac{1}{n}\sum_{i=1}^n Q^i(z_i)\right]\right|\\
					&= \;\;\sup_{f\in \textup{BL}_1(\sZ)}\left| \frac{1}{n}\sum_{i=1}^n\left[\int_{\sZ}f(z_i) \,dP^i(z_i)-\int_{\sZ}f(z_i) \,d Q^i(z_i)\right]\right|\\
					& \stackrel{\eqref{proof:qualitative robustness boottsrap Lemma produktmass und Mischmass}}{= }
					\sup_{f\in \textup{BL}_1(\sZ)}\left| \frac{1}{n}\sum_{i=1}^n\left[\int_{\sZ^{n-1}}\int_{\sZ}f(z_i) \,d  P^i(z_i) \,d\left( \otimes_{j\neq i}P^j(z_j)\right)\right.\right.\\
					&\hspace{15mm}\left.\left.-\int_{\sZ^{n-1}}\int_{\sZ}f(z_i) \,d  Q^i(z_i) \,d\left( \otimes_{j\neq i}Q^j(z_j)\right)\right]\right|\\					
					& = \;\;\sup_{f\in \textup{BL}_1(\sZ)}\left| \frac{1}{n}\sum_{i=1}^n\left[\int_{\sZ^n}f(z_i) \,d \left(\otimes_{j=1}^n P^j(z_j)\right)-\int_{\sZ^n}f(z_i)  \,d \left(\otimes_{j=1}^nQ^j(z_j)\right)\right]\right|\\
					&\leq\;\;\frac{1}{n}\sum_{i=1}^n \sup_{f\in \textup{BL}_1(\sZ)} \left| \int_{\sZ^n}f(z_i)\, d \left(\otimes_{j=1}^n P^j(z_j)\right)-\int_{\sZ^n}f(z_i)  \,d \left(\otimes_{j=1}^nQ^j(z_j)\right)\right|.					
		\end{align*}

Now every function $f\in\textup{BL}_1(\sZ)$ can be identified as a function $\tilde{f}:\sZ^n\rightarrow \sZ$, $(z_1,\ldots,z_n)\mapsto \tilde{f}(z_1,\ldots,z_n):=f(z_i)$. This function is also Lipschitz continuous on $\sZ^n:$ 
\begin{align*}
|\tilde{f}(z_1,\ldots,z_n)-&\tilde{f}(z'_1,\ldots,z'_n)| \;=\;|f(z_i)-f(z'_i)|\\
&\;\leq |f|_1d(z_i,z'_i)\leq |f|_1 (d_{\sZ}(z_1,z'_1)+\ldots+d_{\sZ}(z_i,z'_i)+\ldots+d_{\sZ}(z_n,z'_n)),
\end{align*}
where $d_{\sZ}(z_1,z'_1)+\ldots+d_{\sZ}(z_i,z'_i)+\ldots+d_{\sZ}(z_n,z'_n)$ induces the product topology on $\sZ^n$. That is $\tilde{f}\in\textup{BL}_1(\sZ^n)$. Note that this is also true for every $p$-product metric $d_{n,p}$ in $\sZ^n$, $1\leq p\leq\infty$, as they are strongly equivalent.
Hence, 
\begin{align*}
d_{\textup{BL}}\left( \frac{1}{n}\sum_{i=1}^n P^i, \frac{1}{n}\sum_{i=1}^n Q^i  \right) 
& \leq \frac{1}{n}\sum_{i=1}^n \sup_{g\in \textup{BL}_1(\sZ^n)} \left| \int_{\sZ^n}g \,dP_n-\int_{\sZ^n}g \,d Q_n \right| \\
&\leq
\frac{1}{n}\sum_{i=1}^n \dBL\left(P_n,Q_n\right)\leq \delta,
\end{align*}

which yields the assertion. \hfill$\square$\\


\textbf{Proof of Theorem \ref{THbootstrapindependentnotidentically}:} To prove Theorem \ref{THbootstrapindependentnotidentically} we first use the triangle inequality to split the bounded Lipschitz distance between the distribution of the estimator $S_n$, $n\in\N$, into two parts regarding the distribution of the estimator under the joint distribution $P_n$ of $(Z_1,\ldots, Z_n)$: 	
		\begin{align*}
		\dBL(\sL_{P^*_n}(S_n),\sL_{Q^*_n}(S_n))&\leq \underbrace{\dBL(\sL_{P^*_n}(S_n),\sL_{P_n}(S_n))}_{I}
													+ \underbrace{\dBL(\sL_{P_n}(S_n),\sL_{Q_n^*}(S_n))}_{II}.
		\end{align*}
		
Then the representation of the estimator $S_n$ by the statistical operator $S$ and the continuity of this operator in $P$ together with the Varadarajan property and the independence assumption on the stochastic process yield the assertion.		

First we regard part I: Define the distribution $\PN\in\M(\ZN)$ and let $P^*_{\N}$ be the bootstrap approximation of $\PN$.
Define, for $n\in\N$, the random variables \\
$\W_n: \ZN\rightarrow \sZ^n,\; \W_n=(Z_1,\ldots, Z_n)$, $z_{\N}\mapsto \W_n(z_{\N})=\w_n=(z_1,\ldots, z_n)$, and \\
$\W'_n:\ZN\rightarrow \sZ^n,\;\W'_n=(Z'_1,\ldots, Z'_n)$, $z_{\N}\mapsto \w'_n$, \\
such that $\W_n(\PN) = P_n$ and $\W'_n(\PN^*) = P^*_n$. \\
Denote the bootstrap sample by $\W^*_n:=(Z^*_1,\ldots,Z^*_n)$, $\W^*_n:\ZN\rightarrow \sZ^n$, $z_{\N}\mapsto \w^*_n$. 
\\As Efron's empirical bootstrap is used, the bootstrap sample, which is chosen via resampling with replacement out of $Z_1,\ldots, Z_{\ell}$, $\ell\in\N$, has distribution $Z^*_i\sim \Pn_{\W_{\ell}}=\frac{1}{\ell}\sum_{j=1}^{\ell}\delta_{Z_j}$, $i\in\N$, respectively $\W^*_n:=(Z^*_1,\ldots,Z^*_n)\sim \otimes_{i=1}^n\Pn_{\W_{\ell}}$. The bootstrap approximation of $P_{\ell}$, $\ell\in\N$, is the empirical measure of the bootstrap sample $P^*_{\ell}=\otimes_{i=1}^{\ell} \frac{1}{n}\sum_{j=1}^{n} \delta_{Z^*_j}$.

Further denote the joint distribution of $\W_{\N}$, $\W^*_{\N}$, and $\W'_{\N}$ by $\KN\in\M(\ZN\times\ZN\times\ZN)$. Then, $\KN$ has marginal distributions $\KN(B_1\times \ZN\times\ZN)=\PN(B_1)$ for all $B_1\in \B^{\otimes \N}$, $\KN(\ZN\times B_2\times\ZN)=\otimes_{i\in\N}\Pn_{\W_n}(B_2)$ for all $B_2\in \B^{\otimes \N}$, and $\KN(\ZN\times \ZN\times B_3)=\PN^*(B_3)$ for all $B_3\in \B^{\otimes \N}$.

	Then,
		\begin{equation*}
		\sL_{P_n}(S_n)=S_n(P_n)= S_n\circ \W_n(\PN)\quad \text{and}\quad \sL_{P^*_n}(S_n)=S_n(P^*_n)= S_n \circ \W'_n(\PN^*)
		\end{equation*}
		 and therefore	
		\begin{align*}
		\dBL(\sL_{P^*_n}(S_n),\sL_{P_n}(S_n))		& = \dBL(\sL(S_n\circ W'_n),\sL(S_n\circ W_n)).
		\end{align*}
		
By assumption the coordinate process $(Z_i)_{i\in\N}$ consists of independent random variables, hence we have $P_n=\otimes_{i=1}^{n}P^i$, for $P^i=Z_i(\PN)$, $i\in\N$. 

Moreover $(\sZ,d_{\sZ})$ is assumed to be a totally bounded metric space. Then, due to \citet[Proposition 12]{DudleyZinn1989}, the set $\textup{BL}_1(\sZ,d_{\sZ})$ is a uniform Glivenko-Cantelli class. 
		That is, if $Z_i\sim P$ i.i.d. $i\in\N$, we have for all $\eta>0$: 
		$$\lim_{n\rightarrow\infty}\sup_{P\in\M(\sZ)} \PN \left(\left\{z_{\N}\in\ZN\;|\;\sup_{m\geq n}\dBL(\Pn_{\W_m(z_{\N})},P)>\eta \right\}\right)=0.$$
		 Applying this to the bootstrap sample $(Z^*_1,\ldots, Z^*_m)$, $m\in\N$, which is found by resampling with replacement out of the original sample $(Z_1,\ldots, Z_{n})$, we have, for all $\w_n\in\sZ^n$, 
		\begin{equation*}
		\lim_{n\rightarrow\infty}\sup_{\Pn_{\w_n}\in\M(\sZ)} \otimes_{i\in\N}\Pn_{\w_n} \left(\left\{z_{\N}\in\ZN\;|\;\sup_{m\geq n}\dBL(\Pn_{\W^*_m(z_{\N})},\Pn_{\w_n})>\eta \right\}\right)=0.
		\end{equation*}
		Let $\varepsilon>0$ be arbitrary but fixed. Then, for every $\delta_0 >0$ there is $n_1\in\N$ such that for all $n\geq n_1$ and all $\Pn_{\w_n}\in\M(\sZ)$:
		\begin{equation}\label{proof:bootstrap id Abstandwn*wn}
		\otimes_{i=1}^n\Pn_{\w_n} \left(\left\{\w^*_{n}\in\sZ^n\;|\;\dBL(\Pn_{\w^*_n},\Pn_{\w_n})\leq\frac{\delta_0}{4} \right\}\right)\geq 1-\frac{\varepsilon}{8}.	
		\end{equation}
		
		 And, using the same argumentation for the sequence of random variables $Z'_i$, $i\in\N$, which are i.i.d. and have distribution $\frac{1}{n}\sum_{i=1}^n\delta_{Z^*_i}=\Pn_{\W^*_n}$:	 
		 \begin{equation*}
		\lim_{n\rightarrow\infty}\sup_{\Pn_{\w^*_n}\in\M(\sZ)} P^*_{\N}\left(\left\{z_{\N}\in\ZN\;|\;\sup_{m\geq n}\dBL(\Pn_{\W'_m(z_{\N})},\Pn_{\w^*_n})>\eta \right\}\right)=0.
		\end{equation*}
		 
		 Respectively, for every $\delta_0 >0$ there is $n_2\in\N$ such that for all $n\geq n_2$ and all $\Pn_{\w^*_n}\in\M(\sZ)$:
		\begin{equation}\label{proof:bootstrap id Abstandwn'wn*}
		P^*_n\left(\left\{\w'_{n}\in\sZ^n\;|\;\dBL(\Pn_{\w'_n},\Pn_{\w_n^*})\leq\frac{\delta_0}{2} \right\}\right)\geq 1-\frac{\varepsilon}{8}.	
		\end{equation}
		

		As the process \Z is a strong Varadarajan process by assumption, there exists a probability measure $P\in\M(\sZ)$  
		such that
		\begin{equation*}
		\dBL(\Pn_{\W_n},P)\longrightarrow 0\;\text{almost surely with respect to }\PN,\;n\rightarrow\infty.
		\end{equation*}
		 
		That is, for every $\delta_0 >0$ there is $n_3\in\N$ such that for all $n\geq n_3$:
		\begin{equation}		\label{proof: bootstrap id Varadarajan}	
		P_n\left( \left\{\w_{n}\in\sZ^n\;|\;\dBL(\Pn_{\w_n},P)\leq \frac{\delta_0}{2}\right\}\right)\geq 1-\frac{\varepsilon}{4}.
		\end{equation}		
		
		The continuity of the statistical operator $S:\M(\sZ)\rightarrow H$  in $P\in\M(\sZ)$ yields: for every $\varepsilon >0$ there exists $\delta_0 >0$ such that for all $Q\in\M(\sZ)$: 
		\begin{equation}\label{proof: bootstrap id Stetigkeit Operator}
		\dBL(P,Q)\leq  \delta_0\quad\Rightarrow\quad d_{H}(S(P),S(Q))\leq  \frac{\varepsilon}{4}.
		\end{equation}

		

		As the Prohorov metric $\pi_{d_H}$ is bounded by the Ky Fan metric, see \citet[Theorem 11.3.5]{Dudley1989} we conclude: 
		\begin{align}
		 &\pi_{d_H}(\sL_{P^*_n}(S_n),\sL_{P_n}(S_n))   = \pi_{d_H} (S_n\circ\W'_n,S_n\circ \W_n)\nonumber\\
					&  \leq \inf\left\{\tilde{\varepsilon}>0 \;|\; \KN\left(\left\{d_{H}(S_n\circ \W'_n,S_n\circ \W_n)>\tilde{\varepsilon}\right\}\right)\leq \tilde{\varepsilon}\right\}\nonumber\\	
					&	=  \inf\left\{\tilde{\varepsilon}>0 \;|\;(\W_n,\W^*_n,\W'_n)(\KN)\left(\left\{(\w_n,\w^*_n,\w'_n)\in\sZ^n\times\sZ^n\times\sZ^n\;|\right.\right.\right.\nonumber\\
					&\left.\left.\left.\hspace{4,5cm}\;d_{H}(S_n( \w'_n),S_n(\w_n))>\tilde{\varepsilon},\w^*_n\in\sZ^n\right\}\right)\leq \tilde{\varepsilon}\right\}.	
					\label{proof: bootstrap id Prohorov durch KyFAnMetrik part I}		
		\end{align}

		Due to the definition of the statistical operator $S$, this is equivalent to 	
		\begin{align*}
		  \inf&\left\{\tilde{\varepsilon} >0\;|\;(\W_n,\W^*_n,\W'_n)(\KN)\left(\left\{(\w_n,\w^*_n,\w'_n)\in\sZ^n\times\sZ^n\times\sZ^n\;|\right.\right.\right.\nonumber\\
					&\left.\left.\left.\hspace{4,5cm}
					\;d_{H}(S(\Pn_{\w'_n}),S(\Pn_{\w_n}))>\tilde{\varepsilon}, \w^*_n\in\sZ^n\right\}\right)\leq \tilde{\varepsilon}\right\}.
		\end{align*}

The triangle inequality 
$$d_{H}(S(\Pn_{\w'_n}), S(\Pn_{\w_n}))\leq d_{H}(S(\Pn_{\w'_n}),S(P))+ d_{H}(S(P), S(\Pn_{\w_n})),$$ and the continuity of the statistical operator $S$, see \eqref{proof: bootstrap id Stetigkeit Operator}, then yield, for all $\varepsilon>0$,
		\begin{eqnarray*}
			\lefteqn{(\W_n,\W^*_n,\W'_n)(\KN)\left(\left\{(\w_n,\w^*_n,\w'_n)\in\sZ^n\times\sZ^n\times\sZ^n\;|\;d_{H}(S(\Pn_{\w'_n}),S(\Pn_{\w_n}))>\frac{\varepsilon}{2},\w^*_n\in\sZ^n\right\}\right) }&&\\
						& \leq &  (\W_n,\W^*_n,\W'_n)(\KN)\left(\left\{(\w_n,\w^*_n,\w'_n)\in\sZ^n\times\sZ^n\times\sZ^n\;|\;d_{H}(S(\Pn_{\w'_n}),S(P))> \frac{\varepsilon}{4}\; \right.\right.\\
						&&\hspace{5cm}\left.\left.\text{or}\;
						d_{H}(S(P),S(\Pn_{\w_n}))>\frac{\varepsilon}{4},\w^*_n\in\sZ^n\right\} \right)\\
						& \stackrel{\eqref{proof: bootstrap id Stetigkeit Operator}}{\leq}& (\W_n,\W^*_n,\W'_n)(\KN)\left(\left\{(\w_n,\w^*_n,\w'_n)\in\sZ^n\times\sZ^n\times\sZ^n\;|\; \dBL(\Pn_{\w'_n},P)> \delta_0\; \right.\right.\\
					&&\hspace{5cm}\left.\left.\text{or}\; \dBL(P,\Pn_{\w_n})> \delta_0,\w^*_n\in\sZ^n\right\}\right).
		\end{eqnarray*}
		
		Using the triangle inequality,
		\begin{eqnarray}
		&\quad\dBL(\Pn_{\w'_n},P)\leq \dBL(\Pn_{\w'_n},\Pn_{\w^*_n})+\dBL(\Pn_{\w^*_n},P)\label{proof:qualitative robust bootstrap iid eins triangle a}\\
		\text{and} &\quad \dBL(\Pn_{\w^*_n},P)\leq \dBL(\Pn_{\w^*_n},\Pn_{\w_n})+\dBL(\Pn_{\w_n},P),\label{proof:qualitative robust bootstrap iid eins triangle b}
		\end{eqnarray}
		
		gives for all $n\geq \max\{n_1,n_2,n_3\}$:
		\begin{align*}
		&(\W_n,\W^*_n,\W'_n)(\KN)\left(\left\{(\w_n,\w^*_n,\w'_n)\in\sZ^n\times\sZ^n\times\sZ^n\;|\; \dBL(\Pn_{\w'_n},P)> \delta_0\; \right.\right.\\
					&\hspace{5cm}\left.\left.\text{or}\; \dBL(P,\Pn_{\w_n})> \delta_0,\w^*_n\in\sZ^n\right\}\right)\\
						&\stackrel{\eqref{proof:qualitative robust bootstrap iid eins triangle a}}{\leq} 
							(\W_n,\W^*_n,\W'_n)(\KN)\left(\left\{(\w_n,\w^*_n,\w'_n)\in\sZ^n\times\sZ^n\times\sZ^n\;|\;\dBL(\Pn_{\w'_n},\Pn_{\w^*_n})>\frac{\delta_0}{2}\;\right.\right.\\
						&\left.\left.\hspace{5cm}\text{or}\; \dBL(\Pn_{\w^*_n},P)>  \frac{\delta_0}{2}\;\text{or}\;\dBL(P,\Pn_{\w_n})> \delta_0\right\}\right)\\
						&\stackrel{\eqref{proof:qualitative robust bootstrap iid eins triangle b}}{\leq}
						(\W_n,\W^*_n,\W'_n)(\KN)\left(\left\{(\w_n,\w^*_n,\w'_n)\in\sZ^n\times\sZ^n\times\sZ^n\;|\;\dBL(\Pn_{\w'_n},\Pn_{\w^*_n})>\frac{\delta_0}{2}\;\right.\right.\\
						&\left.\left.\hspace{5cm}\text{or}\; \dBL(\Pn_{\w^*_n},\Pn_{\w_n})>  \frac{\delta_0}{4}\;\text{or}\;\dBL(P,\Pn_{\w_n})> \frac{\delta_0}{4}\right\}\right)\\
						& \leq P^*_{n} \left(\left\{\w'_n\in\sZ^n\;|\; \dBL(\Pn_{\w'_n},\Pn_{\w^*_n})> \frac{\delta_0}{2}\right\}\right) + P_n\left(\left\{\w_n\in\sZ^n\;|\; \dBL(\Pn_{\w_n},P)> \frac{\delta_0}{4}\right\}\right) \\
						&\hspace{4cm}+ \otimes_{i=1}^n\Pn_{\w_n} \left(\left\{\w^*_n\in\sZ^n\;|\; \dBL(\Pn_{\w^*_n},\Pn_{\w_n})> \frac{\delta_0}{4}\right\}\right) \\
						&\stackrel{\eqref{proof:bootstrap id Abstandwn*wn},\eqref{proof:bootstrap id Abstandwn'wn*},\eqref{proof: bootstrap id Varadarajan}}
						{<} \frac{\varepsilon}{8}+\frac{\varepsilon}{4}+ \frac{\varepsilon}{8}\quad =\quad\frac{\varepsilon}{2}.
																	\end{align*}
					

Hence, for all $\varepsilon>0$ there are $n_1,n_2,n_3\in\N$ such that vor all $n\geq \max\{n_1,n_2,n_3\}$, the infimum in \eqref{proof: bootstrap id Prohorov durch KyFAnMetrik part I} is bounded by $\frac{\varepsilon}{2}$. Therefore 
\[
\pi_{d_H} (\sL_{P^*_n}(S_n),\sL_{P_n}(S_n))   < \frac{\varepsilon}{2}.
\]
The equivalence between the Prohorov metric and the bounded Lipschitz metric for Polish spaces, see \citet[Chapter 2, Corollary 4.3]{Huber1981}, yields the existence of  $n_{0,1}\in\N$ such that for all $n\geq n_{0,1}:$
\begin{equation}\label{proof: bootstrap id abschätzung part 1}
\dBL( \sL_{P^*_n}(S_n),\sL_{P_n}(S_n))< \frac{\varepsilon}{2}.
\end{equation}

To prove the convergence of the term in part II, consider the distribution $\QN\in\M(\ZN)$ and let $Q^*_{\N}$ be the bootstrap approximation of $\QN$.
Define, for $n\in\N$, the random variables \\
$\tilde{\W}_n: \ZN \rightarrow \sZ^n,\; \tilde{\W}_n=(\tilde{Z}_1,\ldots, \tilde{Z}_n)$, $z_{\N}\mapsto \tilde{\w}_n$ with distribution  
		$\tilde{\W}_n(\QN)=Q_n$, \\
		$\tilde{\W}'_n:\ZN\rightarrow \sZ^n,\;\tilde{\W}'_n=(\tilde{Z}'_1,\ldots, \tilde{Z}'_n)$, $z_{\N}\mapsto \tilde{\w}'_n$, with distribution 
		$\tilde{\W}'_n(Q^*_{\N})=Q^*_n$, 
		and \\the bootstrap sample
		$\tilde{\W}^*_n:\ZN\rightarrow \sZ^n,\;\tilde{\W}^*_n=(\tilde{Z}^*_1,\ldots, \tilde{Z}^*_n)$, $z_{\N}\mapsto \tilde{\w}^*_n$, with distribution $\otimes_{i=1}^n\Qn_{\tilde{\W}_{\ell}}=\otimes_{i=1}^n\frac{1}{\ell}\sum_{i=1}^{\ell}\delta_{\tilde{Z}_i}.$
		
Moreover let $\KNT\in\M(\ZN\times\ZN\times \ZN\times\ZN)$ denote the joint distribution of $\W_{\N}$, $\tilde{\W}_{\N}$, $\tilde{\W}^*_{\N}$, and $\tilde{\W}'_{\N}$. Then, $\KNT\in\M(\ZN\times\ZN\times \ZN\times\ZN)$ has marginal distributions $\PN $, $\QN$, $\otimes_{i\in\N}\Qn_{\tilde{\W}_n}$, and $ \QN^*$.

First, similar to the argumentation for part I, Efron's bootstrap and \citet[Proposition 12]{DudleyZinn1989} give for $\tilde{\w}_n\in \sZ^n$:
 \begin{equation*}
		 \lim_{n\rightarrow\infty}\sup_{\Qn_{\tilde{\w}_n}\in\M(\sZ)} \otimes_{n\in\N} \Qn_{\tilde{\w}_n}\left(\left\{z_{\N}\in\ZN\;|\;\sup_{m\geq n}\dBL(\Qn_{\tilde{\W}^*_m(z_{\N})},\Qn_{\tilde{\w}_n})>\eta \right\}\right)=0.
		 \end{equation*}
		 
		Hence, for arbitrary, but fixed $\varepsilon >0$, for every $\delta_0 >0$ there is $n_4\in\N$ such that for all $n\geq n_4$ and all $\Qn_{\tilde{\w}_n}\in\M(\sZ)$:
		\begin{equation}\label{proof bootstrap id QsternQn Abstand}
		\otimes_{i=1}^n\Qn_{\tilde{\w}_n}\left(\left\{ \tilde{\w}^*_n\in\sZ^n\;|\;\dBL(\Qn_{\tilde{\w}^*_n},\Qn_{\tilde{\w}_n})\leq \frac{\delta_0}{6} \right\}\right)\geq 1-\frac{\varepsilon}{10}.
		\end{equation}

Further, 
\begin{equation*}
\lim_{n\rightarrow\infty}\sup_{\Qn_{\tilde{\w}^*_n}\in\M(\sZ)} Q^*_{\N}\left(\left\{z_{\N}\in\ZN\;|\;\sup_{m\geq n}\dBL(\Qn_{\tilde{\W}'_m(z_{\N})},\Qn_{\tilde{\w}^*_n})>\eta \right\}\right)=0.	
	\end{equation*}
		
Respectively, for every $\delta_0 >0$ there is $n_5\in\N$ such that for all $n\geq n_5$ and all $\Qn_{\tilde{\w}^*_n}=\frac{1}{n}\sum_{i=1}^n\delta_{\tilde{z}^*_i}\in\M(\sZ)$:
		\begin{equation}\label{proof bootstrap id Q*Qbootstrap Abstand}
		 Q^*_n\left(\left\{ \tilde{\w}'_n\in\sZ^n\;|\;\dBL(\Qn_{\tilde{\w}'_n},\Qn_{\tilde{\w}^*_n})\leq \frac{\delta_0}{6} \right\}\right)\geq 1-\frac{\varepsilon}{10}.
		\end{equation}

	Moreover, as the random variables $Z_i$, $Z_i\sim P^i,$ $i\in\N$, are independent, the bounded Lipschitz distance between the empirical measure and $\frac{1}{n}\sum_{i=1}^n P^i$ can be bounded, due to \citet[Theorem 7]{DudleyZinn1989}.
	As totally bounded spaces are particularly separable, see \citet[below Corollary 1.4.28]{DenkowskiMigorski2003},
	\citet[Proposition 12]{DudleyZinn1989} provides that $\textup{BL}_1(\sZ,d_{\sZ})$ is a uniform Glivenko-Cantelli class. 
	The proof of this proposition does not depend on the distributions of the random variables $Z_i, \;i\in\N$, 
	and is therefore also valid for independent and not necessarily identically distributed random variables. 
	Hence \citet[Theorem 7]{DudleyZinn1989} yields	for all $\eta >0$:
			\[ 
			\lim_{n\rightarrow\infty} \sup_{(P^i)_{i\in\N}\in(\M(\sZ))^{\N}} 
			\PN\left(\left\{z_{\N}\in\ZN\;| \; \sup_{m\geq n} d_{\textup{BL}}\left(\Pn_{\W_m(z_{\N})},\frac{1}{n}\sum_{i=1}^n P^i\right)>\eta\right\}\right)=0,
			\]

	as long as the assumptions of Proposition 12 in \cite{DudleyZinn1989} apply. 
	As $\textup{BL}_1(\sZ,d_{\sZ})$ is bounded, we have 
	$\F_0=\textup{BL}_1(\sZ,d_{\sZ})$, see \citet[page 499, before Proposition 10]{DudleyZinn1989}, 
	hence it is sufficient to show that $\textup{BL}_1(\sZ,d_{\sZ})$ is image admissible Suslin. 
	By assumption $(\sZ,d_{\sZ})$ is totally bounded, 
	hence $\textup{BL}_1(\sZ,d_{\sZ})$ is separable with respect to $\|\cdot\|_{\infty}$, 
	see  \citet[Lemma 3]{StrohrieglHable2016}.
	As $f\in\textup{BL}_1(\sZ,d_{\sZ})$ implies  $\|f\|_{\infty} \leq 1$, the space $\textup{BL}_1(\sZ,d_{\sZ})$ is a bounded subset of $(\C_b(\sZ,d_{\sZ}),\|\cdot\|_{\infty})$, which is due to \citet[Theorem 2.4.9]{Dudley1989} a complete space. Now, $\textup{BL}_1(\sZ,d_{\sZ})$ is a closed subset of $(\C_b(\sZ,d_{\sZ}),\|\cdot\|_{\infty})$ with respect to $\|\cdot\|_{\infty}$. Hence $\textup{BL}_1(\sZ,d_{\sZ})$ is complete, due to \citet[Proposition 1.4.17]{DenkowskiMigorski2003}.
	Therefore $\textup{BL}_1(\sZ,d_{\sZ})$ is separable and complete with respect to $\|\cdot\|_{\infty}$ and 
	particularly a Suslin space, see \citet[p.229]{Dudley1999}. As Lipschitz continuous functions are also equicontinuous, 
	\citet[Theorem 5.28 (c)]{Dudley1999} gives that $\textup{BL}_1(\sZ,d_{\sZ})$ is image admissible Suslin.

	Hence, \citet[Theorem 7]{DudleyZinn1989} yields
		\begin{equation*}\label{dBLPklein}
		\sup_{(P^i)_{i\in\N}\in(\M(\sZ))^{\N}}d_{\textup{BL}}\left(\Pn_{\W_n},\frac{1}{n}\sum_{i=1}^n P^i\right)\longrightarrow 0\;\text{almost surely with respect to }\PN,\;n\rightarrow\infty,
		\end{equation*}
	and
		\begin{equation*}\label{dBLQklein}
		\sup_{(Q^i)_{i\in\N}\in(\M(\sZ))^{\N}}d_{\textup{BL}}\left(\Qn_{\tilde{\W}_n},\frac{1}{n}\sum_{i=1}^n Q^i\right)\longrightarrow 0\;
		\text{almost surely with respect to }\QN,\;n\rightarrow\infty.
		\end{equation*}
		
	That is, there is $n_6\in\N$ such that for all $n \geq n_6$
	\begin{align}
	\quad\quad\; P_n\left(\left\{\w_n\in\sZ^n\;|\;d_{\textup{BL}}\left(\Pn_{\w_n},\frac{1}{n}\sum_{i=1}^n P^i\right)\leq \frac{\delta_0}{6}\right\}\right) \geq 1-\frac{\varepsilon}{10}, 
	\label{proof: bootstrap id Mischmass und empirisches Mass P}\\ 
	\text{and}\quad Q_n\left(\left\{\tilde{\w}_n\in\sZ^n\;|\;d_{\textup{BL}}\left(\Qn_{\tilde{\w}_n},\frac{1}{n}\sum_{i=1}^n Q^i\right)\leq\frac{\delta_0}{6}\right\}\right)\geq 1-\frac{\varepsilon}{10}.
	\label{proof: bootstrap id Mischmass und empirisches Mass Q}
	\end{align}

	Moreover, due to Lemma \ref{Lemma: Zusammenhang Mischmass und Produktmass}, we have
	\begin{equation}
	\dBL(P_n,Q_n)\leq\frac{ \delta_0}{6}\quad\Rightarrow\quad \dBL\left(\frac{1}{n}\sum_{i=1}^n P^i,\frac{1}{n}\sum_{i=1}^n Q^i\right)\leq \frac{\delta_0}{6}.\label{proof: bootstrap id Zusammenhang Mischmass Produktmass}
	\end{equation}
	Then the strong Varadarajan property of \Z 
	yields that 
	there is $n_7\in\N$ such that for all $n\geq n_7:$		
	\begin{equation}\label{proof bootstrap id Varadarajan kleiner n6}
	P_n\left(\left\{\w_n\in\sZ^n\;\vert\;\dBL(\Pn_{\w_n}, P)\leq \frac{\delta_0}{6}\right\}\right) \geq 1-\frac{\varepsilon}{10}.
	\end{equation}
		Similar to the argumentation for part I we conclude, using again the boundedness of the Prohorov metric $\pi_{d_H}$ by the Ky Fan metric, see \citet[Theorem 11.3.5]{Dudley1989}: 
		\begin{align*}
		\lefteqn{\pi_{d_H} (\sL_{P_n}(S_n),\sL_{Q^*_n}(S_n))  \; = \pi_{d_H} (S_n\circ \W_n, S_n\circ\tilde{\W}'_n)}&\\		
					&	=  \inf\{\tilde{\varepsilon}>0 \;|
						\;(\W_n,\tilde{\W}_n,\tilde{\W}^*_n,\tilde{\W}'_n)(\KNT)\left(\left\{(\w_{n},\tilde{\w}_n,\tilde{\w}^*_n,\tilde{\w}'_n)\in \sZ^n\times\sZ^n\times\sZ^n\times\sZ^n\;|\right.\right.\\
						&\hspace{5cm}\left.\left.
						\;d_{H}(S_n(\w_n),S_n(\tilde{\w}'_n))>\tilde{\varepsilon},\;\tilde{\w}_n,\tilde{\w}^*_n\in\sZ^n\right\}\right)\leq \tilde{\varepsilon}\}.
		\end{align*}

		Due to the definition of the statistical operator $S$, this is equivalent to 
		\begin{align*}
		  \inf&\{\tilde{\varepsilon}>0 \;|\;(\W_n,\tilde{\W}_n,\tilde{\W}^*_n,\tilde{\W}'_n)(\KNT)\left(\left\{(\w_{n},\tilde{\w}_n,\tilde{\w}^*_n,\tilde{\w}'_n)\in \sZ^n\times\sZ^n\times\sZ^n\times\sZ^n\;|\right.\right.\\
						&\hspace{5cm}\left.\left.
						\;d_{H}(S(\Pn_{\w_n}),S(\Qn_{\tilde{\w}'_n}))>\tilde{\varepsilon},\;\tilde{\w}_n,\tilde{\w}^*_n\in\sZ^n\right\}\right)\leq \tilde{\varepsilon}\}.
		\end{align*}
		Moreover the triangle inequality yields 		
		\[d_{H}(S(\Pn_{\w_n}), S(\Qn_{\tilde{\w}'_n}))\leq d_{H}(S(\Pn_{\w_n}),S(P))+ d_{H}(S(P), S(\Qn_{\tilde{\w}'_n})) .\]		
		Hence, for all $n\geq \max\{n_4,n_5, n_6,n_7\}$, we obtain		
		\begin{align*}
		&(\W_n,\tilde{\W}_n,\tilde{\W}^*_n,\tilde{\W}'_n)(\KNT)\Big(\Big\{(\w_{n},\tilde{\w}_n,\tilde{\w}^*_n,\tilde{\w}'_n)\in \sZ^n\times\sZ^n\times\sZ^n\times\sZ^n\;|\\
						&\hspace{5cm}\;
					d_{H}(S(\Pn_{\w_n}),S(\Qn_{\tilde{\w}'_n}))> \frac{\varepsilon}{2},\;\tilde{\w}_n,\tilde{\w}^*_n\in\sZ^n\Big\}\Big)&\\
					& \leq (\W_n,\tilde{\W}_n,\tilde{\W}^*_n,\tilde{\W}'_n)(\KNT)\Big(\Big\{(\w_{n},\tilde{\w}_n,\tilde{\w}^*_n,\tilde{\w}'_n)\in \sZ^n\times\sZ^n\times\sZ^n\times\sZ^n\;|\\
						&\hspace{3 cm}\;
						d_{H}(S(\Pn_{\w_n}),S(P))> \frac{\varepsilon}{4}\;\text{or}\;
						d_{H}(S(P),S(\Qn_{\tilde{\w}'_n}))> \frac{\varepsilon}{4},\;\tilde{\w}_n,\tilde{\w}^*_n\in\sZ^n\Big\} \Big).
		\end{align*}

		 The continuity of the statistical operator $S$ in $P$, see 
		\eqref{proof: bootstrap id Stetigkeit Operator}, gives
		\begin{align*}
		&\dBL(P,\Qn_{\tilde{\W}'_n})\leq \delta_0\quad\Rightarrow\quad d_{H}(S(P), S(\Qn_{\tilde{\W}'_n}))\leq\frac{\varepsilon}{4},\\
		\text{and}\quad&\dBL(P,\Pn_{\W_n})\leq \delta_0\quad\Rightarrow\quad d_{H}(S(P), S(\Pn_{\W_n}))\leq\frac{\varepsilon}{4}.
		\end{align*} 
		Further, the triangle inequality yields
		\begin{align}
		\dBL(P,\Qn_{\tilde{\w}'_n}) &\leq 
	\dBL(P,\Pn_{\w_n})+
	\dBL\left(\Pn_{\w_n},\frac{1}{n}\sum_{i=1}^n P^i\right)+
	\dBL\left(\frac{1}{n}\sum_{i=1}^n P^i,\frac{1}{n}\sum_{i=1}^n Q^i\right)\nonumber\\
	&\hspace{5mm}  + 
	\dBL\left(\frac{1}{n}\sum_{i=1}^n Q^i,\Qn_{\tilde{\w}_n}\right)+
	\dBL(\Qn_{\tilde{\w}_n},\Qn_{\tilde{\w}^*_n})+\dBL(\Qn_{\tilde{\w}^*_n},\Qn_{\tilde{\w}'_n}).
	\label{proof: bootstrap id OSternchenP}
		\end{align}
		
		Therefore we conclude, for all $n\geq \max\{n_4,n_5, n_6,n_7\}$,		 
		\begin{align*}
		&(\W_n,\tilde{\W}_n,\tilde{\W}^*_n,\tilde{\W}'_n)(\KNT)\left(\left\{(\w_{n},\tilde{\w}_n,\tilde{\w}^*_n,\tilde{\w}'_n)\in \sZ^n\times\sZ^n\times\sZ^n\times\sZ^n\;|\;
						\right.\right.\\
						&\hspace{3cm}\left.\left.
						d_{H}(S(\Pn_{\w_n}),S(P))> \frac{\varepsilon}{4}\;\text{or}\;d_{H}(S(P),S(\Qn_{\tilde{\w}'_n}))> \frac{\varepsilon}{4},\;\tilde{\w}_n,\tilde{\w}^*_n\in\sZ^n\right\} \right)\\
						& \stackrel{\eqref{proof: bootstrap id Stetigkeit Operator}}{ \leq} (\W_n,\tilde{\W}_n,\tilde{\W}^*_n,\tilde{\W}'_n)(\KNT)\left(\left\{(\w_{n},\tilde{\w}_n,\tilde{\w}^*_n,\tilde{\w}'_n)\in \sZ^n\times\sZ^n\times\sZ^n\times\sZ^n\;|\right.\right.\\
						&\hspace{3cm}\left.\left.\;
					\dBL(\Pn_{\w_n},P)> \delta_0\;\text{or}\;\dBL(P,\Qn_{\tilde{\w}'_n})> \delta_0,\;\tilde{\w}_n,\tilde{\w}^*_n\in\sZ^n \right\}\right)
					\end{align*}
					\begin{align*}	
					& \stackrel{\eqref{proof: bootstrap id OSternchenP}}{\leq}
						(\W_n,\tilde{\W}_n,\tilde{\W}^*_n,\tilde{\W}'_n)(\KNT)\left(\left\{(\w_{n},\tilde{\w}_n,\tilde{\w}^*_n,\tilde{\w}'_n)\in \sZ^n\times\sZ^n\times\sZ^n\times\sZ^n\;|\right.\right.\\
						&\hspace{2cm}\left.\left.
						\;\dBL(\Pn_{\w_n},P)> \frac{\delta_0}{6}\;\text{or}\;	\dBL\left(\Pn_{\w_n},\frac{1}{n}\sum_{i=1}^n P^i\right)> \frac{\delta_0}{6}	\right.\right. \\
   & \hspace{2cm}\left.\left.									
						\;\text{or}\;\dBL\left(\frac{1}{n}\sum_{i=1}^n P^i,\frac{1}{n}\sum_{i=1}^n Q^i\right) > \frac{\delta_0}{6}\;\text{or}\;\dBL\left(\frac{1}{n}\sum_{i=1}^n Q^i,\Qn_{\tilde{\w}_n}\right)> \frac{\delta_0}{6}	\right.\right.\\		
    &\hspace{2cm}\left.\left.															
						\;
						\text{or}\;
							\dBL(\Qn_{\tilde{\w}_n},\Qn_{\tilde{\w}^*_n})> \frac{\delta_0}{6}\;\text{or}\;
							\dBL(\Qn_{\tilde{\w}^*_n},\Qn_{\tilde{\w}'_n})> \frac{\delta_0}{6}
								\right\} \right). 	
																						\end{align*}

%
%
    				
%
%
%
		
		Now, assume $\dBL(P_n,Q_n)\leq \frac{\delta_0}{6}$ , then \eqref{proof: bootstrap id Zusammenhang Mischmass Produktmass} yields $\dBL\left(\frac{1}{n}\sum_{i=1}^n P^i,\frac{1}{n}\sum_{i=1}^n Q^i\right) \leq \frac{\delta_0}{6}$, therefore this term can be omitted. Note that this is only proven for the $p$-product metrics on $\sZ^n$ and not for the metric $d_n$ from \eqref{infmetrik}. For this metric we need a different argumentation, which is stated below the next calculation.
		
		Hence, for all $n\geq\max\{n_4,n_5,n_6,n_7\}$,
		\begin{align*}
		&(\W_n,\tilde{\W}_n,\tilde{\W}^*_n,\tilde{\W}'_n)(\KNT)\left(\left\{(\w_{n},\tilde{\w}_n,\tilde{\w}^*_n,\tilde{\w}'_n)\in \sZ^n\times\sZ^n\times\sZ^n\times\sZ^n\;|\right.\right.\\
		&\hspace{5cm}\left.\left.\;d_{H}(S(\Pn_{\w_n}),S(\Qn_{\tilde{\w}'_n}))> \varepsilon,\;\tilde{\w}_n,\tilde{\w}^*_n\in\sZ^n\right\}\right)\\	
				&\stackrel{\eqref{proof: bootstrap id Zusammenhang Mischmass Produktmass}}{ \leq}
						   (\W_n,\tilde{\W}_n,\tilde{\W}^*_n,\tilde{\W}'_n)(\KNT)\left(\left\{(\w_{n},\tilde{\w}_n,\tilde{\w}^*_n,\tilde{\w}'_n)\in \sZ^n\times\sZ^n\times\sZ^n\times\sZ^n\;|\;
						   		\right.\right.\\
				& \hspace{15mm} \left.\left.
						\dBL(\Pn_{\w_n},P)> \frac{\delta_0}{6}\;\text{or}\;	\dBL\left(\Pn_{\w_n},\frac{1}{n}\sum_{i=1}^n P^i\right)> \frac{\delta_0}{6}\;
						\text{or}\;
							\dBL\left(\frac{1}{n}\sum_{i=1}^n Q^i,\Qn_{\tilde{\w}_n}\right)> \frac{\delta_0}{6}\;
							\right.\right. \\
					& \hspace{15mm}\left.\left.\text{or}\;
							\dBL(\Qn_{\tilde{\w}_n},\Qn_{\tilde{\w}^*_n})> \frac{\delta_0}{6}	\;\text{or}\;
							\dBL(\Qn_{\tilde{\w}^*_n},\Qn_{\tilde{\w}'_n})> \frac{\delta_0}{6}		
								 \right\}\right)\\	
					&  \hspace{3mm}\leq 
						P_n\left(\left\{\w_{n}\in\sZ^n\;|\;\dBL(\Pn_{\w_n},P)> \frac{\delta_0}{6}\right\}\right)\\
						& \hspace{8mm}+
						P_n\left(\left\{\w_{n}\in\sZ^n\;|\;\dBL\left(\Pn_{\w_n},\frac{1}{n}\sum_{i=1}^n P^i\right)> \frac{\delta_0}{6}\right\}\right)\\
				&\hspace{  8mm}
						+Q_n\left(\left\{\tilde{\w}_n\in\sZ^n\;|\;\dBL\left(\frac{1}{n}\sum_{i=1}^n Q^i,\Qn_{\tilde{\w}_n}\right)> \frac{\delta_0}{6}\right\}\right)\\
						&\hspace{8mm}
						+ \otimes_{i=1}^n\Qn_{\tilde{\w}_n}\left(\left\{\tilde{\w}^*_n\in\sZ^n\;|\;\dBL\left(\Qn_{\tilde{\w}_n},\Qn_{\tilde{\w}^*_n}\right)> \frac{\delta_0}{6}\right\}\right)\\
				&\hspace{  8mm}
						+Q^*_n\left(\left\{\tilde{\w}'_n\in\sZ^n\;|\;\dBL\left(\Qn_{\tilde{\w}^*_n},\Qn_{\tilde{\w}'_n}\right)> \frac{\delta_0}{6}\right\}\right)\; 
							\\						 
				& \stackrel{\eqref{proof bootstrap id QsternQn Abstand},
				\eqref{proof bootstrap id Q*Qbootstrap Abstand}
				\eqref{proof: bootstrap id Mischmass und empirisches Mass P},
								\eqref{proof: bootstrap id Mischmass und empirisches Mass Q}, 
									 \eqref{proof bootstrap id Varadarajan kleiner n6}}{ <} 	 
									 \frac{\varepsilon}{10} + \frac{\varepsilon}{10}+ \frac{\varepsilon}{10}+\frac{\varepsilon}{10}+ \frac{\varepsilon}{10}\;=\frac{\varepsilon}{2} .	
									 			 				 						\end{align*}

In order to show the above bound for the metric $d_{n}$, see \eqref{infmetrik}, on $\sZ^n$, we use another variant of the triangle inequality in \eqref{proof: bootstrap id OSternchenP}:
\begin{align}
\dBL(P,\Qn_{\tilde{\w}'_n}) \leq 
	\dBL(P,\Pn_{\w_n})+
	\dBL\left(\Pn_{\w_n},\Qn_{\tilde{\w}_n}\right)+
	\dBL(\Qn_{\tilde{\w}_n},\Qn_{\tilde{\w}^*_n})+\dBL(\Qn_{\tilde{\w}^*_n},\Qn_{\tilde{\w}'_n}).\label{proof:bootstrapid dreieck infmetrik}
\end{align}
Assume $\dBL(P_n,Q_n)\leq \frac{\delta^2_0}{64}$. Then, the strong equivalence between the Prohorov metric and the bounded Lipschitz metric on Polish spaces, see \citet[Chapter 2, Corollary 4.3]{Huber1981}, yields $\pi_{d_n}(P_n,Q_n)\leq \sqrt{\dBL(P_n,Q_n)}\leq \frac{\delta_0}{8}$.   Due to \citet[Theorem 11.6.2]{Dudley1989}, $\pi_{d_n}(P_n,Q_n)\leq \frac{\delta_0}{8}$ implies the existence of a probability measure $\mu\in\M(\sZ^n\times\sZ^n)$ with marginal distributions $P_n$ and $Q_n$, such that $\mu\left(\left\{(\w_n,\tilde{\w}_n)\in\sZ^n\times\sZ^n\;|\;d_n(\w_n,\tilde{\w}_n)>\frac{\delta_0}{8}\right\}\right)\leq \frac{\delta_0}{8}$.
By a simple calculation $d_n(\w_n,\tilde{\w}_n)\leq\frac{\delta_0}{8} $ implies $\pi_{d_n}\left(\frac{1}{n}\sum_{i=1}^n\delta_{z_i},\frac{1}{n}\sum_{i=1}^n\delta_{\tilde{z}_i}\right)\leq \frac{\delta_0}{8}$ and we have:
$$ \mu\left(\left\{(\w_n,\tilde{\w}_n)\in\sZ^n\times\sZ^n\;|\;\pi_{d_n}(\Pn_{\w_n},\Qn_{\tilde{\w}_n})>\frac{\delta_0}{8}\right\}\right)\leq \frac{\delta_0}{8}.$$
Again the equivalence between the metrics $\pi$ and $\dBL$ yields:

$$ \mu\left(\left\{(\w_n,\tilde{\w}_n)\in\sZ^n\times\sZ^n\;|\;\dBL(\Pn_{\w_n},\Qn_{\tilde{\w}_n})>\frac{\delta_0}{4}\right\}\right)\leq \frac{\delta_0}{8}.$$

Now we choose the joint distribution $\tilde{K}_{\N}$ of $\W_{\N}$, $\tilde{\W}_{\N}$, $\tilde{\W}^*_{\N}$, and $\tilde{\W}'_{\N}$ such that the distribution of $(\W_n,\tilde{\W}_n):\ZN\times\ZN\rightarrow\sZ^n\times\sZ^n$ is $\mu\in\M(\sZ^n\times\sZ^n)$. Then we conclude:
		\begin{align*}
		&(\W_n,\tilde{\W}_n,\tilde{\W}^*_n,\tilde{\W}'_n)(\KNT)\left(\left\{(\w_{n},\tilde{\w}_n,\tilde{\w}^*_n,\tilde{\w}'_n)\in \sZ^n\times\sZ^n\times\sZ^n\times\sZ^n\;|\;
						\right.\right.\\
						&\hspace{2cm}\left.\left.
						d_{H}(S(\Pn_{\w_n}),S(P))> \frac{\varepsilon}{4}\;\text{or}\;d_{H}(S(P),S(\Qn_{\tilde{\w}'_n}))> \frac{\varepsilon}{4},\;\tilde{\w}_n,\tilde{\w}^*_n\in\sZ^n\right\} \right)\\	
					& \stackrel{\eqref{proof: bootstrap id Stetigkeit Operator}, \eqref{proof:bootstrapid dreieck infmetrik}}{\leq}
						(\W_n,\tilde{\W}_n,\tilde{\W}^*_n,\tilde{\W}'_n)(\KNT)\left(\left\{(\w_{n},\tilde{\w}_n,\tilde{\w}^*_n,\tilde{\w}'_n)\in \sZ^n\times\sZ^n\times\sZ^n\times\sZ^n\;|\right.\right.\\
						&\hspace{2cm}\left.\left.
						\;\dBL(\Pn_{\w_n},P)> \frac{\delta_0}{4}\;\text{or}\;	\dBL\left(\Pn_{\w_n},\Qn_{\tilde{\w}_n}\right)> \frac{\delta_0}{4}	\right.\right. \\
   & \hspace{2cm}\left.\left.									
						\;\text{or}\;
							\dBL(\Qn_{\tilde{\w}_n},\Qn_{\tilde{\w}^*_n})> \frac{\delta_0}{4}\;\text{or}\;
							\dBL(\Qn_{\tilde{\w}^*_n},\Qn_{\tilde{\w}'_n})> \frac{\delta_0}{4}
								\right\} \right). 	\\
						& \quad \quad\leq \quad\quad
						P_n\left(\left\{\w_{n}\in\sZ^n\;|\;\dBL(\Pn_{\w_n},P)> \frac{\delta_0}{4}\right\}\right)\\
						& \hspace{2cm}+
						\mu\left(\left\{(\w_{n},\tilde{\w}_n)\in\sZ^n\times\sZ^n\;|\;\dBL\left(\Pn_{\w_n},\Qn_{\tilde{\w}_n}\right)> \frac{\delta_0}{4}\right\}\right)\\
						&\hspace{2cm}
						+ \otimes_{i=1}^n\Qn_{\tilde{\w}_n}\left(\left\{\tilde{\w}^*_n\in\sZ^n\;|\;\dBL\left(\Qn_{\tilde{\w}_n},\Qn_{\tilde{\w}^*_n}\right)> \frac{\delta_0}{4}\right\}\right)\\
				&\hspace{ 2cm}
						+Q^*_n\left(\left\{\tilde{\w}'_n\in\sZ^n\;|\;\dBL\left(\Qn_{\tilde{\w}^*_n},\Qn_{\tilde{\w}'_n}\right)> \frac{\delta_0}{4}\right\}\right). 
																						\end{align*}
Now, adapting  the inequalities in \eqref{proof bootstrap id QsternQn Abstand},
\eqref{proof bootstrap id Q*Qbootstrap Abstand}, and
\eqref{proof bootstrap id Varadarajan kleiner n6} in $\varepsilon$ respectively $n$ yields the boundedness of the above term by $\frac{\varepsilon}{2}$ for $\dBL(P_n,Q_n)\leq \frac{\delta^2_0}{64}$ and for all $n\geq\{n_4,n_5,n_7\}$.

 Now we can go on with the proof similar for both kinds of metrics on $\sZ^n$.

The equivalence between the Prohorov metric and the bounded Lipschitz metric on Polish spaces, see \citet[Chapter 2, Corollary 4.3]{Huber1981}, yields the existence of $n_{0,2}\in\N$ such that for all $n\geq n_{0,2}$, $\dBL(P_n,Q_n)\leq\frac{\delta_0}{6}$ (respectively  $\dBL(P_n,Q_n)\leq\frac{\delta^2_0}{64}$) implies
\begin{equation}\label{proof: bootstrap id abschätzung part 2}
\dBL(\sL_{P_n}(S_n),\sL_{Q^*_n}(S_n))<\frac{\varepsilon}{2}.
\end{equation}

Now, \eqref{proof: bootstrap id abschätzung part 1} and \eqref{proof: bootstrap id abschätzung part 2} yield for all $n\geq \max\{n_{0,1},n_{0,2}\}$:
\begin{equation}\label{proof: bootstrap id dBL Bootstrapmasse}
\dBL(\sL_{P^*_n}(S_n),\sL_{Q^*_n}(S_n))<\varepsilon.
\end{equation}

Recall that $\sL_{P_n^*}(S_n)=: \zeta_n$ and $\sL_{Q_n^*}(S_n)=: \xi_n$ are random quantities with values in $\M(H)$. 
	Hence \eqref{proof: bootstrap id dBL Bootstrapmasse} is equivalent to 
	\[ 
	\E \left[\dBL(\sL_{P_n^*}(S_n),\sL_{Q_n^*}(S_n))\right] < \varepsilon, \; \text{for all}\;n\geq \max\{n_{0,1},n_{0,2}\},
	\]
	respectively 
	\[ 
	\E \left[\dBL(\zeta_n,\xi_n)\right] < \varepsilon, \; \text{for all}\;n\geq \max\{n_{0,1},n_{0,2}\}.
	\]
	Therefore, for all $f\in\textup{BL}_1(\M(\sZ))$ and for all $n\geq \max\{n_{0,1},n_{0,2}\}$:
	\begin{align*}
	\left|\int f d(\sL(\zeta_n))-\int f d(\sL(\xi_n)) \right| 
				& = \left|\E f(\zeta_n)-\E f(\xi_n)\right|\;\leq   \E \left| f(\zeta_n)-f(\xi_n)\right|\\
				& \leq  \E \left(\left|f\right|_1 d_{BL} (\zeta_n,\xi_n)\right)	 \;\;\;<  \varepsilon,
	\end{align*}
	
	by a variant of Strassen's Theorem, see \citet[Chapter 2, Theorem 4.2, (2)$\Rightarrow$(1)]{Huber1981}. That is,
	\[
	\dBL(\sL(\sL_{P_n^*}(S_n)),\sL(\sL_{Q_n^*}(S_n)))  < \varepsilon \;\text{for all} \;n\geq \max\{n_{0,1},n_{0,2}\}.
	\]
	Hence for every $\varepsilon >0$ we find $\delta=\frac{\delta_0}{6}$ and $n_0=\max\{n_{0,1},n_{0,2}\}$ such that for all $n\geq n_0$:
		\begin{align*}
	\hspace{3cm} \dBL(P_n,Q_n)<\delta\quad\Rightarrow \quad \dBL(\sL(\sL_{P_n^*}(S_n)),\sL(\sL_{Q_n^*}(S_n))) < \varepsilon,&&
	\end{align*}																										which yields the assertion. \hfill$\square$		
	

	\textbf{Proof of Example 1:}

			Without any restriction we assume $a=0$. Otherwise regard the process $Z_i-a$, $i\in\N$. By assumption, the random variables $Z_i$, $i\in\N$, are independent. 
	Hence $I_B\circ Z_i$, $i\in\N$, are independent, see for example \citet[Theorem 2.10.6]{HoffmannJorgensen1994} 
	for all measurable $B\in\B$, as $I_B$ is a measurable function. 
	According to \citet[Proposition 2.8]{SteinwartHush2009}, \Z satisfies the SLLNE 
	if there is a probability measure $P$ in $\M(\sZ)$ such that 
	$\lim_{n\rightarrow\infty}\frac{1}{n}\sum_{i=1}^n \E_{\mu} I_B\circ Z_i = P(B)$ for all measurable $B\in\B$. 
	Hence:
	\begin{align*}
	\frac{1}{n}\sum_{i=1}^n \E_{\mu} I_B\circ Z_i \; = \; \frac{1}{n}\sum_{i=1}^n \int I_B \,dZ_i(\mu)  \; 
		=  \;\frac{1}{n}\sum_{i=1}^n \int I_Bf_i \,d\lambda^1 &,
	\end{align*}
	where $f_i(x)= \frac{1}{\sqrt{2\pi}}e^{-\frac{1}{2}(x-a_i)^2}$ denotes the density of the 
	normal distribution $N(0,1)$ with respect to the Lebesgue measure $\lambda^1$. Moreover define $g:\R\rightarrow \R$ by
	\begin{align*}
	g(x)=\left\{\begin{array}{ll} 
			e^{-\frac{1}{2}(x+c)^2}, 	& 	x< -c\\
			\frac{1}{\sqrt{2\pi}}					, 		&	-c \leq x\leq c\\
			e^{-\frac{1}{2}(x-c)^2}, 	& 	c< x\\
	\end{array}\right. \quad x\in\R.
	\end{align*}
	
	Therefore $|f_i|\leq |g|$, for all $i\in\N$, $g$ is integrable and due to Lebesgue's Theorem, 
	see for example \citet[Theorem 3.6]{HoffmannJorgensen1994}:	
	\begin{align}
	\lim_{n\rightarrow\infty}\frac{1}{n}\sum_{i=1}^n \int I_Bf_i \,d\lambda^1 \;= \;
					\lim_{n\rightarrow\infty}\int\frac{1}{n}\sum_{i=1}^n  I_B f_i \,d\lambda^1
					\;= \; \int\lim_{n\rightarrow\infty}\frac{1}{n}\sum_{i=1}^n  I_B f_i \,d\lambda^1.\label{noniidexample1}
	\end{align}
	
	We have $f_i\rightarrow f_0$, where $f_0= \frac{1}{\sqrt{2\pi}}e^{-\frac{1}{2}x^2}$ for all $x\in\R$, 
	as $a_i\rightarrow 0$ and therefore  the Lemma of Kronecker, see for example \citet[Theorem 4.9, Equation 4.9.1]{HoffmannJorgensen1994} yields: $\lim_{n\rightarrow\infty}\frac{1}{n}\sum_{i=1}^{n}f_i(x) = f_0(x)$ for all $x\in\X$. 
	
	
	
	
	Now \eqref{noniidexample1} yields the SLLNE:
	
	\[
	\lim_{n\rightarrow\infty}\frac{1}{n}\sum_{i=1}^n \int I_Bf_i \,d\lambda^1 = \int I_B f_0 \,d\lambda^1 = P(B),\;\text{for al}\; B\in\B.\]
																					
	 With  \citet[Zheorem 2]{StrohrieglHable2016} the Varadarajan property is given.						\hfill$\square$

	\textbf{Proof of Example 2:}

	Similar to the proof of Example 1, we first show the SLLNE, that is there exists a probability measure $P\in \M(\sZ)$ such that 
	\[
	\lim_{n\rightarrow\infty}\frac{1}{n}\sum_{i=1}^n \int I_B\circ Z_i\,d\mu= P(B),\;
	\text{for all measurable }B\subset\Omega.
	\]
	Now let $B\subset \Omega$ be an arbitrary measurable set. Then:	
	\begin{align}
	\lim_{n\rightarrow\infty}&\frac{1}{n}\sum_{i=1}^n\int I_B\circ Z_i \,d\mu
				= \lim_{n\rightarrow\infty}\frac{1}{n}\sum_{i=1}^n\int_{\sZ} I_B\,dP^i 
				= \lim_{n\rightarrow\infty}\frac{1}{n}\sum_{i=1}^n\int_{\sZ} I_B\,d 
					[(1-\varepsilon_i)P+\varepsilon_i \tilde{P}^i]\nonumber\\
				&= \lim_{n\rightarrow\infty}\frac{1}{n}\sum_{i=1}^n\int_{\sZ} I_B\,dP -
					\lim_{n\rightarrow\infty}\frac{1}{n}\sum_{i=1}^n\varepsilon_i\int_{\sZ} I_B\,d P+ 
					\lim_{n\rightarrow\infty}\frac{1}{n}\sum_{i=1}^n\varepsilon_i\int_{\sZ} I_B\,d\tilde{P}^i.
					\label{shrinkinge}
	\end{align}
	
	As,  $0\leq \frac{1}{n}\sum_{i=1}^n \varepsilon_i \int I_B\,dP \leq \frac{1}{n}\sum_{i=1}^n \varepsilon_i$ 
	and $\varepsilon_i\rightarrow 0$, we have 
	\[
	\lim_{n\rightarrow\infty}\frac{1}{n}\sum_{i=1}^n \varepsilon_i\int I_B\,dP\leq 
	\lim_{n\rightarrow\infty}\frac{1}{n}\sum_{i=1}^n \varepsilon_i \longrightarrow 0,\; n\rightarrow\infty
	\]
	and similarly
	\[
	\lim_{n\rightarrow\infty}\frac{1}{n}\sum_{i=1}^n \varepsilon_i\int I_B\,d\tilde{P}^i\leq 
	\lim_{n\rightarrow\infty}\frac{1}{n}\sum_{i=1}^n \varepsilon_i \longrightarrow 0,\quad n\rightarrow\infty.
	\]
	Hence \eqref{shrinkinge} yields
	\begin{align*}
	\lim_{n\rightarrow\infty}\frac{1}{n}\sum_{i=1}^n I_B\circ Z_i  = 
	\lim_{n\rightarrow\infty}\frac{1}{n}\sum_{i=1}^n \int I_B\,dP 	= P(B)
	\end{align*}
	and therefore, due to \citet[Theorem 2]{StrohrieglHable2016}, the assertion.				\hfill$\square$

\textbf{Proof of Corollary \ref{COBootstrapidSVM}:}

	Due to Example 2, the stochastic process is a Varadarajan process.
	\citet[Theorem 3.2]{HableChristmann2011} ensures the continuity of the statistical operator $S:\M(\sZ)\rightarrow H,\; P\mapsto f_{L^*,P,\lambda}$ for a fixed value $\lambda\in(0,\infty)$. Moreover \citet[Corollary 3.4]{HableChristmann2011} yields the continuity of the estimator 
	$S_n: \sZ^n\rightarrow H,
	\; D_n\mapsto f_{L^*,D_n,\lambda}$ for every fixed $\lambda\in(0,\infty)$. Hence for fixed $\lambda >0$ the bootstrap approximation of the SVM estimator is qualitatively robust, for the given assumptions. 
	  Moreover the proof of Theorem \ref{THbootstrapindependentnotidentically}, equation \eqref{proof: bootstrap id dBL Bootstrapmasse}, and the equivalence between between bounded Lipschitz metric and Prokhorov distance yield: for every $\varepsilon >0$ there is $\delta >0$ such that there is $n_0\in\N$ such that for all $n\geq n_0$ and if $\dBL(P_n,Q_n)\leq\delta$:
\begin{equation*}
\pi(\sL_{P^*_n}(S_n),\sL_{Q^*_n}(S_n))<\varepsilon\;\text{almost surely}.
\end{equation*}
Similarly to the proof of the qualitative robustness in \citet[Theorem 4]{StrohrieglHable2016} we get: for every $\varepsilon >0$ there is $n_{\varepsilon}$, such that for all $n\geq n_{\varepsilon}$: 
 \[
 \|f_{L^*,D_n,\lambda_n}-f_{L^*,D_n,\lambda_0}\|_H\leq\frac{\varepsilon}{3}.\] 
 And the same argumentation as in the proof of the qualitative robustness of the SVM estimator for the non-i.i.d. case in \citet[Theorem 4]{StrohrieglHable2016} for the cases $n_{0}\leq n\leq n_{\varepsilon}$ and $n>n_{\varepsilon}$ yields the assertion. \hfill $\square$

\subsection{Proofs of Section \ref{subsec:bootstrapbetamixing}}
		\textbf{Proof of Theorem \ref{THbootstrapqualitativerobustnessmixing}:}
	
		\textbf{Proof of Theorem \ref{THbootstrapqualitativerobustnessmixing}:}	Let $P^*_{\N},Q^*_{\N}\in\M(\ZN)$ be the bootstrap approximations of the true distribution $\PN$ and the contaminated distribution $Q_{\N}$. First, the triangle inequality yields:
		\begin{align*}
		\lefteqn{\dBL(\sL_{P^*_n}(S_n),\sL_{Q^*_n}(S_n))}&\nonumber\\ 
				&\leq \underbrace{\dBL(\sL_{P^*_n}(S_n),\sL_{P_n}(S_n))}_{I}
													+ \underbrace{\dBL(\sL_{P_n}(S_n),\sL_{Q_n}(S_n))}_{II}
													+ \underbrace{\dBL(\sL_{Q_n}(S_n),\sL_{Q^*_n}(S_n))}_{III}.
		\end{align*}

		First, we regard the term in part II. 
		Let $\sigma(Z_i),\;i\in\N$, be the $\sigma$-algebra generated by $Z_i$. 
		Due to the assumptions on the mixing process 
		$\sum_{m>n}\alpha (\sigma(Z_1,\ldots,Z_i),\sigma(Z_{i+m},\ldots),\PN) = \SO(n^{-\gamma}),\; i\in\N, \gamma >0$, the sequence 
		$\left(\alpha (\sigma(Z_1,\ldots,Z_i),\sigma(Z_{i+m},\ldots),\mu)\right)_{m\in\N}$ 
		is a null sequence. Moreover it is bounded by the definition of the $\alpha$-mixing coefficient which, due to the strong stationarity, does not depend on $i$. 
		Therefore
			\begin{align*}
			\frac{1}{n^2}\sum_{i=1}^n\sum_{j=1}^n \alpha ((Z_i)_{i\in\N},\PN,i,j) & =\frac{1}{n^2}\sum_{i=1}^n\sum_{j=1}^n \alpha (\sigma(Z_i),\sigma(Z_j),\PN) \\
					&\leq \frac{2}{n^2}\sum_{i=1}^n\sum_{j\geq i}^n \alpha (\sigma(Z_i),\sigma(Z_j),\PN) \\
					&\leq \frac{2}{n^2}\sum_{i=1}^n\sum_{j\geq i}^n \alpha (\sigma(Z_1,\ldots,Z_i),\sigma(Z_j,\ldots),\PN) \\
					& = \frac{2}{n^2}\sum_{i=1}^n\sum_{\ell=0}^{n-i} \alpha (\sigma(Z_1,\ldots,Z_i),\sigma(Z_{i+\ell},\ldots),\PN) \\
					& \stackrel{stationarity}{\leq} \frac{2}{n}\sum_{\ell =0}^{n}\alpha (\sigma(Z_1,\ldots,Z_i),\sigma(Z_{i+\ell},\ldots),\PN),\;i\in\N\\
					&\longrightarrow 0,\;n\rightarrow\infty.
			\end{align*}
			Hence, the process is weakly $\alpha$-$bi$-mixing with respect to $\PN$, see Definition \ref{Def: bootstrap weakly alpha bi mixing}. Due to the stationarity assumption, the process \Z is additionally  asymptotically mean stationary, that is $\lim_{n\rightarrow\infty}\frac{1}{n}\sum_{i=1}^{n} \E I_B\circ Z_i =P(B)$ for all $B\in\A$ for a probability measure $P$. Therefore the process satisfies the WLLNE, see \citet[Proposition 3.2]{SteinwartHush2009}, and therefore is a weak Varadarajan process, see \citet[Theorem 2]{StrohrieglHable2016}.
			
		As the process is assumed to be a Varadarajan process and due to the assumptions on the sequence of estimators $(S_n)_{n\in\N},$ qualitative robustness of $\left(S_n\right)_{n\in\N}$ is ensured by \citet[Theorem 1]{StrohrieglHable2016}. Together with the equivalence between the Prohorov metric and the bounded Lipschitz metric for Polish spaces, see \citet[Chapter 2, Corollary 4.3]{Huber1981}, it follows:
		
		For every $\varepsilon >0$ there is $\delta>0$ such that for all $n\in\N$ and for all $Q_n\in\M(\sZ^n)$ we have:		
		\begin{equation*}
		\dBL(P_n,Q_n)< \delta \quad \Rightarrow \quad \dBL(\sL_{P_n}(S_n),\sL_{Q_n}(S_n)) < \frac{\varepsilon}{3}.
		\end{equation*}
		 This implies
		 \begin{equation} \label{proof:bootstrap beta mixing part II}
		 \E \left[\dBL(\sL_{P_n}(S_n),\sL_{Q_n}(S_n)) \right]< \frac{\varepsilon}{3}.
		 \end{equation}
		Hence the convergence of the term in part II is shown.

		
To prove the convergence of the term in part I, consider the distribution $\PN\in\M(\ZN)$ and let $P^*_{\N}$ be the bootstrap approximation of $\PN$, via the blockwise bootstrap.
Define, for $n\in\N$, the random variables \\
$\W_n: \ZN\rightarrow \sZ^n,\; \W_n=(Z_1,\ldots, Z_n)$, $z_{\N}\mapsto \w_n=(z_1,\ldots,z_n)$, and \\
$\W'_n:\ZN\rightarrow \sZ^n,\;\W'_n=(Z'_1,\ldots, Z'_n)$, $z_{\N}\mapsto \w'_n$, \\
such that $\W_n(\PN) = P_n$ and $\W'_n(\PN^*) = P^*_n$. \\
Moreover denote the bootstrap sample by $\W^*_n:\ZN\rightarrow \sZ^n$, $\W^*_n:=(Z^*_1,\ldots,Z^*_n)$, $z_{\N}\mapsto \w^*_n$, and the distribution of $\W^*_n$ by  $\overline{P}_n$. 
The blockwise bootstrap approximation of $P_{m}$, $m\in\N$, is $P^*_{m}=\otimes_{j=1}^{m} \frac{1}{n}\sum_{i=1}^{n} \delta_{Z^*_i}$, $m\in\N$. Note that the sample $Z^*_1,\ldots,Z^*_n$ depends and on the blocklength $b(n)$ and on the number of blocks $\ell(n)$.

Further denote the joint distribution of $\W_{\N}$, $\W^*_{\N}$, and $\W'_{\N}$ by $\KN\in\M(\ZN\times\ZN\times\ZN)$. Then, $\KN$ has marginal distributions $\KN(B_1\times \ZN\times\ZN)=\PN(B_1)$ for all $B_1\in \B^{\otimes \N}$, $\KN(\ZN\times B_2\times\ZN)=\overline{P}_{\N}(B_2)$ for all $B_2\in \B^{\otimes \N}$, and $\KN(\ZN\times \ZN\times B_3)=\PN^*(B_3)$ for all $B_3\in \B^{\otimes \N}$.\\	
		Then,
		\begin{equation*}
		\sL_{P_n}(S_n)=S_n(P_n)= S_n\circ \W_n(\PN)\quad \text{and}\quad \sL_{P^*_n}(S_n)=S_n(P^*_n)= S_n \circ \W'_n(\PN^*)
		\end{equation*}
		 and therefore	
		\begin{align*}
		\dBL(\sL_{P^*_n}(S_n),\sL_{P_n}(S_n))		& = \dBL(\sL(S_n\circ W'_n),\sL(S_n\circ W_n)).
		\end{align*}		
		By assumption we have $0\leq z_i\leq 1$, $i\in\N$. Hence $Z_i(z_{\N})=z_i\in [0,1]$, i.\,e. $\sZ=[0,1]$, which is a totally bounded metric space. 		
		Therefore the set $\textup{BL}_1([0,1])$ is a uniform Glivenko-Cantelli class, 
		due to \citet[Proposition 12]{DudleyZinn1989}. 
		Similar to part I of the proof of Theorem \ref{THbootstrapindependentnotidentically}, the blockwise bootstrap structure and the Glivenko-Cantelli property yield:
		
		 \begin{equation*}
		\lim_{n\rightarrow\infty}\sup_{\Pn_{\w^*_n}\in\M(\sZ)} P^*_{\N}\left(\left\{z_{\N}\in\ZN\;|\;\sup_{m\geq n}\dBL(\Pn_{\W'_m(z_{\N})},\Pn_{\w^*_n})>\eta \right\}\right)=0.
		\end{equation*}
		 
		 Respectively, for fixed $\varepsilon >0$, for every $\delta_0 >0$ there is $n_1\in\N$ such that for all $n\geq n_1$ and all $\Pn_{\w^*_n}\in\M(\sZ)$:
		\begin{equation}\label{proof:bootstrap alphamixing R Abstandwn'wn*}
		P^*_n\left(\left\{\w'_{n}\in\sZ^n\;|\;\dBL(\Pn_{\w'_n},\Pn_{\w^*_n})\leq\frac{\delta_0}{2} \right\}\right)\geq 1-\frac{\varepsilon}{6}.	
		\end{equation}

		Regard the process 
		$G_n(t)=\frac{1}{\sqrt{n}}\sum_{i=1}^n I_{\{Z^*_i\leq t\}} - \frac{1}{\sqrt{n}}\sum_{i=1}^nI_{\{Z_i\leq t\}}$, 
		$t\in\R$.
		Due to the assumptions on the process and on the moving block bootstrap, Theorem 2.3 in \cite{Peligrad1998} 
		yields the almost sure convergence in distribution to a Brownian bridge $G$: 
		\begin{equation}
		\frac{1}{\sqrt{n}}\sum_{i=1}^n I_{\{Z^*_i\leq t\}} - \frac{1}{\sqrt{n}}\sum_{i=1}^nI_{\{Z_i\leq t\}} \longrightarrow_{\D} G(t),\quad t\in\R
		\end{equation}
		 almost surely with respect to $\PN$, $n\rightarrow\infty$,  in the Skorohod topology on $D[0,1]$.
		 Here $\longrightarrow_{\D}$ indicates convergence in distribution and $D[0,1]$ denotes the space of cadlag functions on $[0,1]$, for details see for example \citet[p. 121]{Billingsley2013convergence}.
		
		This is equivalent to 		 
		 \begin{equation*}
		\frac{1}{\sqrt{n}}\sum_{i=1}^n I_{\{Z^*_i\leq t\}} - \frac{1}{\sqrt{n}}\sum_{i=1}^nI_{\{Z_i\leq t\}} \longrightarrow_{\D} G(t),\;
		\text{almost surely with respect to } \PN,\;n\rightarrow\infty,
		\end{equation*}
		for all continuity points $t$ of $G$, see \citet[(12.14), p. 124]{Billingsley2013convergence}. 
		
		Multiplying by $\frac{1}{\sqrt{n}}$ yields for any fixed continuity point $t\in\R:$		
		\begin{equation*}
		\frac{1}{n}\sum_{i=1}^n I_{\{Z^*_i\leq t\}} - \frac{1}{n}\sum_{i=1}^nI_{\{Z_i\leq t\}} -\frac{1}{\sqrt{n}}G(t)\longrightarrow_{\D} 0\;\text{almost surely with respect to } \PN,\;n\rightarrow\infty.
		\end{equation*}
		
		As convergence in distribution to a finite constant implies convergence in probability, see for example \citet[Theorem 2.7(iii)]{VanderVaart1988}, and as $\frac{1}{\sqrt{n}}G(t)\rightarrow 0$ in probability, for all $t\in\R$:	
		\begin{equation*}
		\frac{1}{n}\sum_{i=1}^n I_{\{Z^*_i\leq t\}} - \frac{1}{n}\sum_{i=1}^nI_{\{Z_i\leq t\}} \longrightarrow_{P} 0 \; \text{almost surely with respect to } \PN,\;n\rightarrow\infty,
		\end{equation*}
for all continuity points $t$ of $G$, where $\longrightarrow_{P}$ denotes the convergence in probability.

		Hence, \citet[Theorem 11.12]{Dudley1989} yields the convergence of the corresponding probability measures:		
		\begin{equation*}
		\dBL \left(\frac{1}{n}\sum_{i=1}^n \delta_{Z^*_i}, \frac{1}{n}\sum_{i=1}^n \delta_{Z_i}\right)\longrightarrow _P 0 
		\;\text{almost surely with respect to } \PN,\;n\rightarrow\infty.
		\end{equation*}
		
		 Respectively		 
		 \begin{equation*}
		\dBL (\Pn_{\W^*_n}, \Pn_{\W_n})\longrightarrow _P 0 \;\text{almost surely with respect to } \PN,\;n\rightarrow\infty.
		\end{equation*}
		
		 Define the set $B_n= \left\{\w_n\in\sZ^n\;|\;\dBL (\Pn_{\W^*_n}, \Pn_{\w_n})\longrightarrow_P 0,\; n\rightarrow\infty\right\}$. Hence,
		 \begin{equation}\label{proof: bootstrap alphamixing BnP}
		 P_n(B_n)=\PN\left(\left\{z_{\N}\in\ZN\;|\;\W_n(z_{\N})\in B_n\right\}\right)=1
\end{equation}		  
		and, for all $\w_n\in B_n$, there is $n_{2,\w_n}\in\N$ such that for all $n\geq n_{2,\w_n}\in\N$:
		\begin{equation}\label{proof: bootstrap mixingR Konvergenz PwStern Pw 1}
		\overline{P}_n\left(\left\{\w^*_{n}\in\sZ^n\;|\;\dBL\left(\Pn_{\w^*_n}, \Pn_{\w_n}\right)>\frac{\delta_0}{4}\right\}\right) < \frac{\varepsilon}{6}.
		\end{equation}
		
		By assumption we have $0\leq z_i\leq 1$, $i\in\N$. 
		Hence the space of probability measures $\left\{\Pn_{\w_n}\;|\; \w_n\in [0,1]^n\right\}$ is a subset of $\M([0,1])$ and therefore tight, as [0,1] is a compact space, see e.\,g. \cite[Example 13.28]{Klenke}. 
		Then Prohorov's Theorem, 
		see for example \citet[Theorem 5.1]{Billingsley2013convergence} yields relative compactness of $\M([0,1],\dBL)$ 
		and in particular the relative compactness of the set $\left\{\Pn_{\w_n}\;|\;\w_n\in [0,1]^n\right\}$. As $\M([0,1],\dBL)$ is a complete space, see \citet[Theorem 11.5.5]{Dudley1989}, relative compactness equals total boundedness.
		That is, there exists a finite dense subset $\tilde{\mathcal{P}}$ of $\left\{\Pn_{\w_n}\;|\; \w_n\in [0,1]^n\right\}$ 
		such that for all $\rho>0$ and $\Pn_{\w_n}\in\left\{\Pn_{\w_n}\;|\; \w_n\in [0,1]^n\right\}$ 
		there is $\tilde{P}_{\rho}\in \tilde{\mathcal{P}}$  such that		
		\begin{align}\label{proof:  bootstrap mixingR M(P) totally bounded}
		\dBL(\tilde{P}_{\rho}, \Pn_{\w_n})\leq\rho.
		\end{align}
		
		The triangle inequality yields:		
		\[\dBL\left(\Pn_{\w^*_n}, \Pn_{\w_n}\right)\leq 
		\dBL\left(\Pn_{\w^*_n}, \tilde{P}_{\rho}\right)+\dBL\left(\tilde{P}^{\rho}, \Pn_{\w_n}\right).\]
		
		Define $\rho=\frac{\delta_0}{4}$. Then \eqref{proof: bootstrap mixingR Konvergenz PwStern Pw 1} yields for every $\tilde{P}_{\rho}\in\tilde{\mathcal{P}}$ the existence of an integer $n\geq n_{2,\tilde{P}}\in\N$ such that, for all $n\geq n_{2,\tilde{P}}$ and all $\w_n\in B_n$:	
		\begin{align*}
		\overline{P}_n & \left(\left\{\w^*_{n}\in\sZ^n\;|\;\dBL\left(\Pn_{\w^*_n}, \Pn_{\w_n}\right)>\frac{\delta_0}{2}\right\}\right)\\
		 & \;\;\leq \;\;\overline{P}_n\left(\left\{\w^*_{n}\in\sZ^n\;|\;\dBL\left(\Pn_{\w^*_n}, \tilde{P}_{\rho}\right)>\frac{\delta_0}{4}\;\text{or}\;\dBL\left(\tilde{P}_{\rho}, \Pn_{\w_n}\right)>\frac{\delta_0}{4}\right\}\right)\\
		 & \stackrel{\eqref{proof:  bootstrap mixingR M(P) totally bounded}}{ \leq} \overline{P}_n\left(\left\{\w^*_{n}\in\sZ^n\;|\;\dBL\left(\Pn_{\w^*_n}, \tilde{P}_{\rho}\right)>\frac{\delta_0}{4}\right\}\right)
		 \;\stackrel{\eqref{proof: bootstrap mixingR Konvergenz PwStern Pw 1}}{< }\;\frac{\varepsilon}{6}. 
		\end{align*}
		
		Hence, for all $n\geq n_2:=\max_{\tilde{P}\in\mathcal{\tilde{P}}}\{ n_{2,\tilde{P}}\}$ and for all $\w_n\in B_n$, we have:		
		\begin{align}\label{proof: bootstrap mixingR Konvergenz PwStern Pw}
		\sup_{\Pn_{\w_n}\in\M(\sZ)} \overline{P}_n&\left(\left\{\w^*_{n}\in\sZ^n\;|\;\dBL\left(\Pn_{\w^*_n}, \Pn_{\w_n}\right)>\frac{\delta_0}{2}\right\}\right)<\frac{\varepsilon}{6}.
		\end{align}		
		
		 Due to the uniform continuity of the operator $S$, for every $\varepsilon >0$ there is $\delta_0 >0$ such that for all $P,Q\in\M(\sZ)$ :
		 \begin{equation}\label{proof:bootstrap alpha mmixing R uniform continuit part II}
		 \dBL(P,Q)\leq\delta_0\quad\Rightarrow \quad d_{H}(S(P),S(Q))\leq \frac{\varepsilon}{3}.
\end{equation}		  
		 

		Moreover, the triangle inequality yields:
		\begin{equation}\label{proof:bootstrap alphamixing R triangle w'wn}
		\dBL(\Pn_{\w'_n},\Pn_{\w_n})\leq \dBL(\Pn_{\w'_n},\Pn_{\w^*_n})+\dBL(\Pn_{\w^*_n},\Pn_{\w_n}).
\end{equation}

		Again we use the relation between the Prohorov metric $\pi_{d_H}$ and the Ky Fan metric, \citet[Theorem 11.3.5]{Dudley1989}: 
		\begin{align*}
		\pi_{d_H} \left(\sL_{P^*_n}(S_n),\right.&\left.\sL_{P_n}(S_n)\right)  \; = \pi_{d_H}(S_n\circ\W'_n,S_n\circ \W_n)\\
					&  \leq \inf\left\{\tilde{\varepsilon}>0 \;|\; \KN\Big(\Big\{			
					\;d_{H}(S_n\circ \W'_n,S_n\circ \W_n)>\tilde{\varepsilon},\w^*_{\N}\in\ZN\Big\}\Big)\leq \tilde{\varepsilon}\right\}\\	
					&	=  \inf\left\{\tilde{\varepsilon}>0 \;|\;(\W_n,\W^*_n,\W'_n)(\KN)\left(\left\{(\w_n,\w^*_n,\w'_n)\in\sZ^n\times\sZ^n\times\sZ^n\;|\right.\right.\right.\\
					&\hspace{25mm}\left.\left.\left.
					\;d_{H}(S_n( \w'_n),S_n(\w_n))>\tilde{\varepsilon},\w^*_n\in\sZ^n\right\}\right)\leq \tilde{\varepsilon}\right\}.	
		\end{align*}

		Due to the definition of the statistical operator $S$, this is equivalent to 	
		\begin{align*}
		  \inf&\{\tilde{\varepsilon}>0 \;|\;(\W_n,\W^*_n,\W'_n)(\KN)\left(\left\{(\w_n,\w^*_n,\w'_n)\in\sZ^n\times\sZ^n\times\sZ^n\;|
		  \right.\right.\\
					&\hspace{45mm}\left.\left.
					\;d_{H}(S(\Pn_{\w'_n}),S(\Pn_{\w_n})>\tilde{\varepsilon},\w^*_n\in\sZ^n\right\}\right)\leq \tilde{\varepsilon}\}.
		\end{align*}

		Due to the uniform continuity of $S$, see \eqref{proof:bootstrap alpha mmixing R uniform continuit part II}, we obtain, for all $n\geq \max\{n_1,n_2\}:$
		\begin{eqnarray*}
			\lefteqn{(\W_n,\W^*_n,\W'_n)(\KN)\left(\left\{(\w_n,\w^*_n,\w'_n)\in\sZ^n\times\sZ^n\times\sZ^n\;|\;d_{H}(S(\Pn_{\w'_n}),S(\Pn_{\w_n}))>\frac{\varepsilon}{3},\w^*_n\in\sZ^n\right\}\right) }&&\\
				&\stackrel{\eqref{proof:bootstrap alpha mmixing R uniform continuit part II}}{ \leq}& (\W_n,\W^*_n,\W'_n)(\KN)\left(\left\{(\w_n,\w^*_n,\w'_n)\in\sZ^n\times\sZ^n\times\sZ^n\;|\;\dBL(\Pn_{\w'_n},\Pn_{\w_n})> \delta_0,\w^*_n\in\sZ^n\right\}\right)\\
				& = &(\W_n,\W^*_n,\W'_n)(\KN)\left(\left\{(\w_n,\w^*_n,\w'_n)\in\sZ^n\times\sZ^n\times\sZ^n\;|\right.\right.\\
				&& \hspace{10mm} \left.\left.
				\;\{\w_n\notin B_n,\;\dBL(\Pn_{\w'_n},\Pn_{\w_n})> \delta_0\}\;\text{or}\;\{\w_n\in B_n,\;\dBL(\Pn_{\w'_n},\Pn_{\w_n})> \delta_0\},\w^*_n\in\sZ^n\right\}\right)\\				
				& \leq &(\W_n,\W^*_n,\W'_n)(\KN)\left(\left\{(\w_n,\w^*_n,\w'_n)\in\sZ^n\times\sZ^n\times\sZ^n\;|\right.\right.\\
				&& \hspace{25mm} \left.\left.
				\;\w_n\notin B_n,\;\dBL(\Pn_{\w'_n},\Pn_{\w_n})> \delta_0,\w^*_n\in\sZ^n\right\}\right)\\
				&&\hspace{5mm}+ 
					(\W_n,\W^*_n,\W'_n)(\KN)\left(\left\{(\w_n,\w^*_n,\w'_n)\in\sZ^n\times\sZ^n\times\sZ^n\;|\right.\right.\\
				&& \hspace{25mm} \left.\left.
				\;\w_n\in B_n,\;\dBL(\Pn_{\w'_n},\Pn_{\w_n})> \delta_0,\w^*_n\in\sZ^n\right\}\right)\\
				&\stackrel{\eqref{proof: bootstrap alphamixing BnP}}{=} &
					(\W_n,\W^*_n,\W'_n)(\KN)\left(\left\{(\w_n,\w^*_n,\w'_n)\in\sZ^n\times\sZ^n\times\sZ^n\;|\right.\right.\\
				&& \hspace{25mm} \left.\left.
				\;\w_n\in B_n,\;\dBL(\Pn_{\w'_n},\Pn_{\w_n})> \delta_0,\w^*_n\in\sZ^n\right\}\right).
				\end{eqnarray*}
				
				The triangle inequality, \eqref{proof:bootstrap alphamixing R triangle w'wn}, then yields for all $n\geq \max\{n_1,n_2\}$:
				\begin{align*}
					(\W_n,\W^*_n,\W'_n)&(\KN) \left(\left\{(\w_n,\w^*_n,\w'_n)\in\sZ^n\times\sZ^n\times\sZ^n|\;
				\w_n\in B_n,\dBL(\Pn_{\w'_n},\Pn_{\w_n})> \delta_0,\w^*_n\in\sZ^n\right\}\right)\\					
				&\stackrel{\eqref{proof:bootstrap alphamixing R triangle w'wn}}{\leq} 
				(\W_n,\W^*_n,\W'_n)(\KN)\Big(\Big\{(\w_n,\w^*_n,\w'_n)\in\sZ^n\times\sZ^n\times\sZ^n\;|\;\{\w_n\in B_n\;\Big.\Big.\\
				&\Big.\Big.\hspace{20mm}
				\text{and}\;\dBL(\Pn_{\w'_n},\Pn_{\w^*_n})> \frac{\delta_0}{2}\} \;\text{or}\;\{\w_n\in B_n \;\text{and}\;\dBL(\Pn_{\w^*_n},\Pn_{\w_n})> \frac{\delta_0}{2}\}  \Big\}\Big)\\	
				& \;\;\leq 	P^*_{n}\left(\left\{\w'_n\in\sZ^n\;|\;\w_n\in\B_n,\;\dBL(\Pn_{\w'_n},\Pn_{\w^*_n})>\frac{\delta_0}{2}\right\}\right)\\
				&\hspace{10mm}+ \overline{P}_{n}\left(\left\{\w^*_n\in\sZ^n\;|\;\w_n\in\B_n,\;\dBL(\Pn_{\w^*_n},\Pn_{\w_n})>\frac{\delta_0}{2}\right\}\right)\\
				&\stackrel{ \eqref{proof:bootstrap alphamixing R Abstandwn'wn*},\eqref{proof: bootstrap mixingR Konvergenz PwStern Pw 1} }{<}	\frac{\varepsilon}{6}+\frac{\varepsilon}{6} \quad = \quad\frac{\varepsilon}{3}.		
		\end{align*}

	The equivalence between the Prohorov metric and the bounded Lipschitz metric on Polish spaces, see \citet[Chapter 2, Corollary 4.3]{Huber1981}, yields the existence of $\tilde{n}_1$ such that for every 
	$n\geq \tilde{n}_1:$
		\begin{equation*}\label{proof: bootstrap mixingR abschätzung part 1}
		\dBL( \sL_{P^*_n}(S_n),\sL_{P_n}(S_n))< \frac{
		\varepsilon}{3}.
		\end{equation*}

And therefore
		\begin{equation}
			 \label{proof:bootstrap beta mixing convergence part I}
			\E \left[\dBL\left(\sL_{P_n^*}(S_n), \sL_{P_n}(S_n)\right) \right]<\frac{
			\varepsilon}{3} .	  
		\end{equation}

		For the convergence of the term in part III the same argumentation as for part I can be applied, as the assumptions on $\QN$ and $Q^*_{\N}$ are the same as for $\PN$ and $\Pn^*_{\N}$. In particular for every $\varepsilon >0$ there is $\tilde{n}_2\in\N$ such that for all $n\geq \tilde{n}_2$: 
		\begin{equation*}	
		\dBL\left(\sL_{Q^*_n}(S_n), \sL_{Q_n}(S_n)\right) <  \frac{
		\varepsilon}{3}, 
		\end{equation*}
		respectively		
		\begin{equation}\label{proof:bootstrap beta mixing convergence part III}
		\E \left[\dBL\left(\sL_{Q^*_n}(S_n), \sL_{Q_n}(S_n)\right) \right]<  \frac{
		\varepsilon}{3} .
		\end{equation}
		
		Hence, \eqref{proof:bootstrap beta mixing part II}, \eqref{proof:bootstrap beta mixing convergence part I}, and \eqref{proof:bootstrap beta mixing convergence part III} yield, for all $n\geq \max\{\tilde{n}_1,\tilde{n}_2\}$:
		
		\begin{equation*}
		\E \left[\dBL\left(\sL_{P^*_n}(S_n), \sL_{Q_n^*}(S_n)\right)\right]< \frac{\varepsilon}{3}+\frac{\varepsilon}{3}+\frac{\varepsilon}{3}   = \varepsilon.
		\end{equation*}
		
		As $\sL_{P^*_n}(S_n)$ and $\sL_{Q^*_n}(S_n)$ are random variables itself we have, due to \citet[Chapter 2 Theorem 4.2, (2)$\Rightarrow$(1)]{Huber1981}, for all $n\geq\max\{\tilde{n}_1,\tilde{n}_2\}$:
		\begin{equation*}
		\dBL\left(\sL(\sL_{P^*_n}(S_n))	, \sL(\sL_{Q^*_n}(S_n))\right)<
		\varepsilon.
		\end{equation*}		
		
		Hence, for all $\varepsilon>0$ there is $\delta >0$ such that there is $n_0=\max\{\tilde{n}_1,\tilde{n}_2\}\in\N$ such that, for all $n\geq n_0$:
	 \begin{equation*}
 	d_{\textup{BL}}(P_n,Q_n)< \delta \;\Rightarrow\; d_{\textup{BL}}(\sL(\sL_{P^*_n}(S_n)),\sL(\sL_{Q^*_n}(S_n))) < \varepsilon
 	\end{equation*}
	 
	 and therefore the assertion.  \hfill$\square$


\textbf{Proof of Theorem \ref{THbootstraphigherdimensions}:}
																								
\textbf{Proof of Theorem \ref{THbootstraphigherdimensions}:} 
The proof follows the same lines as the proof of Theorem \ref{THbootstrapqualitativerobustnessmixing} and therefore we only state the different steps. Again we start with the triangle inequality:
\begin{align*}
		\lefteqn{\dBL(\sL_{P^*_n}(S_n),\sL_{Q^*_n}(S_n))}&\nonumber\\ 
				&\leq \underbrace{\dBL(\sL_{P^*_n}(S_n),\sL_{P_n}(S_n))}_{I}
													+ \underbrace{\dBL(\sL_{P_n}(S_n),\sL_{Q_n}(S_n))}_{II}
													+ \underbrace{\dBL(\sL_{Q_n}(S_n),\sL_{Q^*_n}(S_n))}_{III}.
		\end{align*}

To proof the convergence of the term in part II, we need the weak Varadarajan property of the stochastic process. Due to the definition $\alpha(\sigma(Z_1,\ldots,Z_i),\sigma(Z_{i+\ell},\ldots),\mu) \leq 2$ for all $\ell\in\N$, $i\in\N$, and obviously:
\begin{equation}\label{proof: bootstrapmixingRdAbschaetzungalpha}
\alpha(\sigma(Z_1,\ldots,Z_i),\sigma(Z_{i+\ell},\ldots),\PN) \leq \ell+1,\; \ell>0.
\end{equation}

Hence, due to the strong stationarity of the stochastic process, we have:
\begin{align*}
				\lefteqn{\frac{1}{n^2}\sum_{i=1}^n\sum_{j=1}^n \alpha ((Z_i)_{i\in\N},\PN,i,j) \;=\;\frac{1}{n^2}\sum_{i=1}^n\sum_{j=1}^n \alpha (\sigma(Z_i),\sigma(Z_j),\PN) 
					}& \\
					&\quad \leq\quad  \frac{2}{n^2}\sum_{i=1}^n\sum_{j\geq i}^n \alpha (\sigma(Z_i),\sigma(Z_j),\PN)\\
					&\quad	\leq \quad\frac{2}{n^2}\sum_{i=1}^n\sum_{j\geq i}^n \alpha (\sigma(Z_1,\ldots,Z_i),\sigma(Z_j,\ldots),\PN) \\
					& \quad= \quad \frac{2}{n^2}\sum_{i=1}^n\sum_{\ell=0}^{n-i} \alpha (\sigma(Z_1,\ldots,Z_i),\sigma(Z_{i+\ell},\ldots),\PN) \\
					& \stackrel{\textup{stationarity}}{\leq} \frac{2}{n}\sum_{\ell=0}^{n}\alpha (\sigma(Z_1,\ldots,Z_i),\sigma(Z_{i+\ell},\ldots),\PN),\; i\in\N\\
					&\quad =\quad \frac{2}{n}\sum_{\ell=0}^{n}\left(\alpha (\sigma(Z_1,\ldots,Z_i),\sigma(Z_{i+\ell},\ldots),\PN)\right)^{\frac{1}{2}}\left(\alpha (\sigma(Z_1,\ldots,Z_i),\sigma(Z_{i+\ell},\ldots),\PN)\right)^{\frac{1}{2}},\; i\in\N\\
					& \quad\stackrel{\eqref{proof: bootstrapmixingRdAbschaetzungalpha}}{\leq} \quad
					\frac{2}{n}\sum_{\ell=0}^{n}(\ell+1)\left(\alpha (\sigma(Z_1,\ldots,Z_i),\sigma(Z_{i+\ell},\ldots),\PN)\right)^{\frac{1}{2}},\; i \in\N\\
					&\quad\stackrel{\eqref{proof: bootstrapmixing Rd Bedingung an Pozess}}{\longrightarrow} 0, \;n\rightarrow\infty.
			\end{align*}
			
Now, the same argumentation as in the proof of Theorem \ref{THbootstrapqualitativerobustnessmixing} yields the weak Varadarajan property and therefore, for all $\varepsilon>0$,			
 	
		
		 \begin{equation} \label{proof:bootstrap beta mixing part II higher dimensions}
		 \E \left[\dBL(\sL_{P_n}(S_n),\sL_{Q_n}(S_n))\right] < \frac{\varepsilon}{3}.
		 \end{equation}

Regarding the term in part I, we use a central limit theorem for the blockwise bootstrapped empirical process by \citet[Corollary 1 and remark] {Buehlmann1994} to show its convergence.  	
Again, regard the distribution $\PN\in\M(\ZN)$ and let $P^*_{\N}$ be the bootstrap approximation of $\PN$, via the blockwise bootstrap.
Define, for all $n\in\N$, the random variables \\
$\W_n: \ZN\rightarrow \sZ^n,\; \W_n=(Z_1,\ldots, Z_n)$, $z_{\N}\mapsto \w_n$, and \\
$\W'_n:\ZN\rightarrow \sZ^n,\;\W'_n=(Z'_1,\ldots, Z'_n)$, $z_{\N}\mapsto \w'_n$, \\
such that $\W_n(\PN) = P_n$ and $\W'_n(\PN^*) = P^*_n$. \\
Moreover denote the bootstrap sample by $\W^*_n:\ZN\rightarrow\sZ^n$, $\W^*_n:=(Z^*_1,\ldots,Z^*_n)$, $z_{\N}\mapsto \w^*_n$, and the distribution of $\W^*_n$ by  $\overline{P}_n$. The bootstrap approximation of $P_{m}$ is $P^*_{m}=\otimes_{j=1}^{m} \frac{1}{n}\sum_{i=1}^{n} \delta_{Z^*_i}=\otimes_{j=1}^{m}\Pn_{\W^*_n}$, $m\in\N$, by definition of the bootstrap procedure. Note that the sample $Z^*_1,\ldots,Z^*_n$ depends and on the blocklength $b(n)$ and on the number of blocks $\ell(n)$.

Further denote the joint distribution of $\W_{\N}$, $\W^*_{\N}$, and $\W'_{\N}$ by $\KN\in\M(\ZN\times\ZN\times\ZN)$. Then, $\KN$ has marginal distributions $\KN(B_1\times \ZN\times\ZN)=\PN(B_1)$ for all $B_1\in \B^{\otimes \N}$, $\KN(\ZN\times B_2\times\ZN)=\overline{P}_{\N}(B_2)$ for all $B_2\in \B^{\otimes \N}$, and $\KN(\ZN\times \ZN\times B_3)=\PN^*(B_3)$ for all $B_3\in \B^{\otimes \N}$.
				
		Then,
		\begin{equation*}
		\sL_{P_n}(S_n)=S_n(P_n)= S_n\circ \W_n(\PN)\quad \text{and}\quad \sL_{P^*_n}(S_n)=S_n(P^*_n)= S_n \circ \W'_n(\PN^*)
		\end{equation*}
		 and therefore	
		\begin{align*}
		\dBL(\sL_{P^*_n}(S_n),\sL_{P_n}(S_n))		& = \dBL(\sL(S_n\circ W'_n),\sL(S_n\circ W_n)).
		\end{align*}		
		As $\sZ=[0,1]^d$ is compact, it is in particular totally bounded. 		
		Hence the set $\textup{BL}_1(\sZ,d_{\sZ})$ is a uniform Glivenko-Cantelli class, 
		due to \citet[Proposition 12]{DudleyZinn1989}. 
		Similar to part I of the proof of Theorem \ref{THbootstrapqualitativerobustnessmixing}, the bootstrap structure and the Glivenko-Cantelli property given above yield for arbitrary, but fixed $\varepsilon>0$:\\
for every $\delta_0 >0$ there is $n_0\in\N$ such that, for all $n\geq n_0$ and all $\Pn_{\w^*_n}\in\M(\sZ)$,
		\begin{equation*}\label{proof:bootstrap alphamixing Rd Abstandwn'wn*}
		P^*_n\left(\left\{\w'_{n}\in\sZ^n\;|\;\dBL(\Pn_{\w'_n},\Pn_{\w^*_n})\leq\frac{\delta_0}{2} \right\}\right)\geq 1-\frac{\varepsilon}{6}.	
		\end{equation*}

 	Now, regard the empirical process of $(Z_1,\ldots, Z_n)$. Set 
 	$\ttt=(t_1,\ldots, t_d)\in\R^d$. 
 	Moreover $\ttt<\mathbf{b}$ means $t_i<b_i$ for all $i\in\{1,\ldots,d\}$. 
 	Hence we can define the empirical process and the blockwise bootstrapped empirical process by
 	
 	$$\frac{1}{n}\sum_{i=1}^nI_{\{Z_i\leq\ttt\}}\quad\text{and}\quad\frac{1}{n}\sum_{i=1}^nI_{\{Z^*_i\leq \ttt\}}.$$
 	
 	Regard the process 
 	$G_n(\ttt)=\frac{1}{\sqrt{n}}\sum_{i=1}^n I_{\{Z^*_i\leq \ttt\}} -\frac{1}{\sqrt{n}}\sum_{i=1}^n I_{\{Z_i\leq \ttt\}}$, $\ttt\in [0,1]^d$.
		Now, due to the assumptions on the stochastic process and on the moving block bootstrap, \citet[Corollary 1 and remark]{Buehlmann1994} 
		yields the almost sure convergence in distribution to a Gaussian process $G$: 
		\begin{equation*}
		\frac{1}{\sqrt{n}}\sum_{i=1}^n I_{\{Z^*_i\leq \ttt\}} -\frac{1}{\sqrt{n}}\sum_{i=1}^n I_{\{Z_i\leq \ttt\}} 
		\longrightarrow_{\D} G(\ttt),\quad \ttt\in [0,1]^d,
		\end{equation*}	
		 almost surely with respect to $\PN$, $n\rightarrow\infty$, in the (extended) Skorohod topology on $D^d([0,1])$.\\
		 The space $D^d([0,1])$ is a generalization of the space of cadlag functions on $[0,1]$, see \citet[Chapter 12]{Billingsley2013convergence}, and consists of functions $f:[0,1]^d\rightarrow\R$. A detailed description of this space and the extended Skorohod topology can be found in \cite{Straf1972, Straf1969} and \cite{BickelWichura1971}. The definition of the space $D^{d}([0,1])$ can, for example, be found in \citet[Chapter 3]{BickelWichura1971}.
		 
		\citet[Lemma 5.4]{Straf1972} yields, that the above convergence in the Skorohod topology is equivalent to the convergence for all continuity points $\ttt$ of $G$. Hence,		 
		 \begin{equation*}
		\frac{1}{\sqrt{n}}\sum_{i=1}^n I_{\{Z^*_i\leq \ttt\}} -\frac{1}{\sqrt{n}}\sum_{i=1}^n I_{\{Z_i\leq \ttt\}} 
		\longrightarrow_{\D} G(\ttt) \;\text{almost surely with respect to }\PN,\;n\rightarrow\infty,
		\end{equation*}
		for all continuity points $\ttt$ of $G$. 
		
		Multiplying by $\frac{1}{\sqrt{n}}$ yields, for every continuity point $\ttt$ of $G$,	
		\begin{equation*}
		\frac{1}{n}\sum_{i=1}^n I_{\{Z^*_i\leq \ttt\}} - \frac{1}{n}\sum_{i=1}^nI_{\{Z_i\leq \ttt\}} -\frac{1}{\sqrt{n}}G(\ttt)
		\longrightarrow_{\D} 0\;\text{almost surely with respect to }\PN,\;n\rightarrow\infty.
		\end{equation*}

As convergence in distribution to a constant implies convergence in probability, see e.\,g. \citet[Theorem 2.7(iii)]{VanderVaart1988} and as $\frac{1}{\sqrt{n}}G(\ttt)$ converges in probability to $0$, for all fixed continuity points $\ttt\in [0,1]^d$ of $G$:
		\begin{equation*}
		\frac{1}{n}\sum_{i=1}^n I_{\{Z^*_i\leq \ttt\}} - \frac{1}{n}\sum_{i=1}^nI_{\{Z_i\leq \ttt\}} \longrightarrow_{P} 0  \;\text{almost surely with respect to }\PN,\;n\rightarrow\infty.
		\end{equation*}		
		This yields the convergence of the corresponding probability measures, see for example 
		\citet[Chapter 29]{Billingsley2008probability} for a theory on $\R^d$:		
		\begin{equation*}
		\dBL (\frac{1}{n}\sum_{i=1}^n \delta_{Z^*_i}, \frac{1}{n}\sum_{i=1}^n \delta_{Z_i})\longrightarrow _P 0 \;\text{almost surely with respect to }\PN,\;n\rightarrow\infty,
		\end{equation*}	
			 
		 respectively		 
		 \begin{equation*}
		\dBL (\Pn_{\W^*_n}, \Pn_{\W_n})\longrightarrow _P 0 \;\text{almost surely with respect to }\PN,\;n\rightarrow\infty.
		\end{equation*}	
			 
		 As the space $[0,1]^d$ is compact, we can use an argumentation similar to the proof of Theorem \ref{THbootstrapqualitativerobustnessmixing}. Then, for every $\varepsilon >0$, there is $n_1\in\N$ such that for all $n\geq n_1$
		\begin{equation*}
		\dBL\left(\sL_{P^*_n}(S_n), \sL_{P_n}(S_n)\right) <\frac{
		\varepsilon}{3},
		\end{equation*}		
		respectively,		
		\begin{equation}\label{proof:bootstrap beta mixing convergence part I higher dimensions}
		\E \left[\dBL\left(\sL_{P^*_n}(S_n), \sL_{P_n}(S_n)\right) \right]<\frac{
		\varepsilon}{3} .
		\end{equation}			
		The convergence of the term in part III follows simultaneously to part I for the distributions $Q_{\N}$ and $Q^*_{\N}$. Hence,	
		for every $\varepsilon >0$, there is $n_2\in\N$ such that for all $n\geq n_2$		
		\begin{equation}\label{proof:bootstrap beta mixing convergence part III higher dimensions}
		\E \left[\dBL\left(\sL_{Q^*_n}(S_n), \sL_{Q_n}(S_n)\right) \right]<  \frac{
		\varepsilon}{3} .
		\end{equation}		
		The combination of \eqref{proof:bootstrap beta mixing part II higher dimensions}, \eqref{proof:bootstrap beta mixing convergence part I higher dimensions}, and \eqref{proof:bootstrap beta mixing convergence part III higher dimensions} yields for all $n\geq \max\{n_1,n_2\}$:	
		\begin{equation*}
		\E \left[\dBL\left(\sL_{P^*_n}(S_n), \sL_{Q^*}(S_n)\right)\right]< \frac{\varepsilon}{3}+\frac{\varepsilon}{3}+\frac{\varepsilon}{3}   = \varepsilon.
		\end{equation*}
		
		As $\sL_{P^*_n}(S_n)$ and $\sL_{Q^*_n}(S_n)$ are random variables itself we have, due to \citet[Chapter 2, Theorem 4.2, (2)$\Rightarrow$(1)]{Huber1981}, for all $n\geq\max\{n_1,n_2\}:$
		\begin{equation*}
		\dBL\left(\sL(\sL_{P^*_n}(S_n))	, \sL(\sL_{Q^*_n}(S_n))\right)<
		\varepsilon.
		\end{equation*}		
		
		Hence, for all $\varepsilon>0$ there is $\delta >0$ such that there is $n_0=\max\{n_1,n_2\}\in\N$ such that for all $n\geq n_0:$
	 \begin{equation*}
 	d_{\textup{BL}}(P_n,Q_n)< \delta \;\Rightarrow\; d_{\textup{BL}}(\sL(\sL_{P^*_n}(S_n)),\sL(\sL_{Q^*_n}(S_n))) < \varepsilon.
 	\end{equation*}
	 
	 This yields the assertion.  \hfill$\square$

\bibliographystyle{abbrvnat}
\bibliography{literature}

\begin{thebibliography}{54}
\providecommand{\natexlab}[1]{#1}
\providecommand{\url}[1]{\texttt{#1}}
\expandafter\ifx\csname urlstyle\endcsname\relax
  \providecommand{\doi}[1]{doi: #1}\else
  \providecommand{\doi}{doi: \begingroup \urlstyle{rm}\Url}\fi

\bibitem[Beutner and Z\"ahle(2016)]{BeutnerZaehle2016}
E.~Beutner and H.~Z\"ahle.
\newblock Functional delta-method for the bootstrap of quasi-{H}adamard
  differentiable functionals.
\newblock \emph{Electron. J. Stat.}, 10, 2016.

\bibitem[Bickel and Wichura(1971)]{BickelWichura1971}
P.~J. Bickel and M.~J. Wichura.
\newblock Convergence criteria for multiparameter stochastic processes and some
  applications.
\newblock \emph{Ann. Math. Statist.}, 42:\penalty0 1656--1670, 1971.

\bibitem[Billingsley(1995)]{Billingsley2008probability}
P.~Billingsley.
\newblock \emph{Probability and measure}.
\newblock Wiley Series in Probability and Mathematical Statistics. John Wiley
  \& Sons, Inc., New York, third edition, 1995.

\bibitem[Billingsley(1999)]{Billingsley2013convergence}
P.~Billingsley.
\newblock \emph{Convergence of probability measures}.
\newblock Wiley Series in Probability and Statistics: Probability and
  Statistics. John Wiley \& Sons, Inc., New York, second edition, 1999.

\bibitem[Boente et~al.(1982)Boente, Fraiman, and Yohai]{BoenteFairmanYohai1982}
G.~Boente, R.~Fraiman, and V.~J. Yohai.
\newblock Qualitative robustness for general stochastic processes.
\newblock Technical report, Department of Statistics, University of Washington,
  1982.

\bibitem[Boente et~al.(1987)Boente, Fraiman, and Yohai]{BoenteFairmanYohai1987}
G.~Boente, R.~Fraiman, and V.~J. Yohai.
\newblock Qualitative robustness for stochastic processes.
\newblock \emph{The Annals of Statistics}, 15\penalty0 (3):\penalty0
  1293--1312, 1987.

\bibitem[Bradley(2005)]{Bradley2005}
R.~C. Bradley.
\newblock Basic properties of strong mixing conditions. {A} survey and some
  open questions.
\newblock \emph{Probab. Surv.}, 2:\penalty0 107--144, 2005.

\bibitem[Bradley(2007{\natexlab{a}})]{Bradley20071}
R.~C. Bradley.
\newblock \emph{Introduction to strong mixing conditions. {V}ol. 1}.
\newblock Kendrick Press, Heber City, UT, 2007{\natexlab{a}}.

\bibitem[Bradley(2007{\natexlab{b}})]{Bradley20072}
R.~C. Bradley.
\newblock \emph{Introduction to strong mixing conditions. {V}ol. 2}.
\newblock Kendrick Press, Heber City, UT, 2007{\natexlab{b}}.

\bibitem[Bradley(2007{\natexlab{c}})]{Bradley20073}
R.~C. Bradley.
\newblock \emph{Introduction to strong mixing conditions. {V}ol. 3}.
\newblock Kendrick Press, Heber City, UT, 2007{\natexlab{c}}.

\bibitem[B{\"u}hlmann(1994)]{Buehlmann1994}
P.~B{\"u}hlmann.
\newblock Blockwise bootstrapped empirical process for stationary sequences.
\newblock \emph{Ann. Statist.}, 22\penalty0 (2):\penalty0 995--1012, 1994.

\bibitem[B{\"u}hlmann(1995)]{Buehlmann1995}
P.~B{\"u}hlmann.
\newblock The blockwise bootstrap for general empirical processes of stationary
  sequences.
\newblock \emph{Stochastic Process. Appl.}, 58\penalty0 (2):\penalty0 247--265,
  1995.

\bibitem[Bustos(1980)]{Bustos1980}
O.~Bustos.
\newblock On qualitative robustness for general processes.
\newblock unpublished manuscript, 1980.

\bibitem[Christmann et~al.(2011)Christmann, Salibian-Barrera, and
  Van~Aelst]{Christmann2011}
A.~Christmann, M.~Salibian-Barrera, and S.~Van~Aelst.
\newblock On the stability of bootstrap estimators.
\newblock \emph{arXiv preprint arXiv:1111.1876}, 2011.

\bibitem[Christmann et~al.(2013)Christmann, Salibi{\'a}n-Barrera, and
  Van~Aelst]{ChristmannVanAelst2013}
A.~Christmann, M.~Salibi{\'a}n-Barrera, and S.~Van~Aelst.
\newblock Qualitative robustness of bootstrap approximations for kernel based
  methods.
\newblock In \emph{Robustness and complex data structures}, pages 263--278.
  Springer, Heidelberg, 2013.

\bibitem[Cox(1981)]{Cox1981}
D.~D. Cox.
\newblock metrics on stochastic processes and qualitative robustness.
\newblock Technical report, Department of Statistics, University of Washington,
  1981.

\bibitem[Cuevas and Romo(1993)]{CuevasRomo1993}
A.~Cuevas and J.~Romo.
\newblock On robustness properties of bootstrap approximations.
\newblock \emph{J. Statist. Plann. Inference}, 37\penalty0 (2):\penalty0
  181--191, 1993.

\bibitem[Denkowski et~al.(2003)Denkowski, Mig{\'o}rski, and
  Papageorgiou]{DenkowskiMigorski2003}
Z.~Denkowski, S.~Mig{\'o}rski, and N.~S. Papageorgiou.
\newblock \emph{An introduction to nonlinear analysis: applications}.
\newblock Kluwer Academic Publishers, Boston, MA, 2003.

\bibitem[Doukhan(1994)]{Doukhan1994}
P.~Doukhan.
\newblock \emph{Mixing}.
\newblock Springer, New York, 1994.

\bibitem[Dudley(1989)]{Dudley1989}
R.~M. Dudley.
\newblock \emph{Real Analysis and Probability}.
\newblock Chapman$\&$Hall, New York, 1989.

\bibitem[Dudley(2014)]{Dudley1999}
R.~M. Dudley.
\newblock \emph{Uniform central limit theorems}, volume~63 of \emph{Cambridge
  Studies in Advanced Mathematics}.
\newblock Cambridge University Press, Cambridge, 2014.

\bibitem[Dudley et~al.(1991)Dudley, Gin{\'e}, and Zinn]{DudleyZinn1989}
R.~M. Dudley, E.~Gin{\'e}, and J.~Zinn.
\newblock Uniform and universal {G}livenko-{C}antelli classes.
\newblock \emph{J. Theoret. Probab.}, 4\penalty0 (3):\penalty0 485--510, 1991.

\bibitem[Efron(1979)]{Efron1979}
B.~Efron.
\newblock Bootstrap methods: another look at the jackknife.
\newblock \emph{Ann. Statist.}, 7\penalty0 (1):\penalty0 1--26, 1979.

\bibitem[Efron and Tibshirani(1993)]{Efron1993}
B.~Efron and R.~J. Tibshirani.
\newblock \emph{An introduction to the bootstrap}, volume~57 of
  \emph{Monographs on Statistics and Applied Probability}.
\newblock Chapman and Hall, New York, 1993.

\bibitem[Hable and Christmann(2011)]{HableChristmann2011}
R.~Hable and A.~Christmann.
\newblock On qualitative robustness of support vector machines.
\newblock \emph{Journal of Multivariate Analysis}, 102:\penalty0 993--1007,
  2011.

\bibitem[Hampel(1968)]{Hampel1968}
F.~R. Hampel.
\newblock \emph{Contributions to the theory of robust estimation}.
\newblock PhD thesis, Univ. California, Berkeley, 1968.

\bibitem[Hampel(1971)]{Hampel1971}
F.~R. Hampel.
\newblock A general qualitative definition of robustness.
\newblock \emph{Annals of Mathematical Statistics}, 42:\penalty0 1887--1896,
  1971.

\bibitem[Hoffmann-J{\o}rgensen(1994)]{HoffmannJorgensen1994}
J.~Hoffmann-J{\o}rgensen.
\newblock \emph{Probability with a view toward statistics. {V}ol. {I}}.
\newblock Chapman \& Hall Probability Series. Chapman \& Hall, New York, 1994.

\bibitem[Huber(1981)]{Huber1981}
P.~J. Huber.
\newblock \emph{Robust statistics}.
\newblock John Wiley \& Sons Inc., New York, 1981.

\bibitem[Jure{\v{c}}kov{\'a} and Picek(2006)]{JureckovaPicek2006}
J.~Jure{\v{c}}kov{\'a} and J.~Picek.
\newblock \emph{Robust statistical methods with {$R$}}.
\newblock Chapman \& Hall/CRC, Boca Raton, FL, 2006.

\bibitem[Klenke(2013)]{Klenke}
A.~Klenke.
\newblock \emph{Probability theory: a comprehensive course}.
\newblock Springer Science \& Business Media, 2013.

\bibitem[Kr\"atschmer et~al.(2017)Kr\"atschmer, Schied, and
  Z\"ahle]{Zaehle2017}
V.~Kr\"atschmer, A.~Schied, and H.~Z\"ahle.
\newblock Domains of weak continuity of statistical functionals with a view
  toward robust statistics.
\newblock \emph{J. Multivariate Anal.}, 158:\penalty0 1--19, 2017.

\bibitem[K{\"u}nsch(1989)]{Kuensch1989}
H.~R. K{\"u}nsch.
\newblock The jackknife and the bootstrap for general stationary observations.
\newblock \emph{Ann. Statist.}, 17\penalty0 (3):\penalty0 1217--1241, 1989.

\bibitem[Lahiri(2003)]{Lahiri2003}
S.~N. Lahiri.
\newblock \emph{Resampling methods for dependent data}.
\newblock Springer Series in Statistics. Springer, New York, 2003.

\bibitem[Liu and Singh(1992)]{LiuSingh1992}
R.~Y. Liu and K.~Singh.
\newblock Moving blocks jackknife and bootstrap capture weak dependence.
\newblock In \emph{Exploring the limits of bootstrap ({E}ast {L}ansing, {MI},
  1990)}, Wiley Ser. Probab. Math. Statist. Probab. Math. Statist., pages
  225--248. Wiley, New York, 1992.

\bibitem[Maronna et~al.(2006)Maronna, Martin, and Yohai]{MaronnaMartin2006}
R.~A. Maronna, R.~D. Martin, and V.~J. Yohai.
\newblock \emph{Robust statistics}.
\newblock Wiley Series in Probability and Statistics. John Wiley \& Sons Ltd.,
  Chichester, 2006.

\bibitem[Naik-Nimbalkar and Rajarshi(1994)]{Naik1994}
U.~V. Naik-Nimbalkar and M.~B. Rajarshi.
\newblock Validity of blockwise bootstrap for empirical processes with
  stationary observations.
\newblock \emph{Ann. Statist.}, 22\penalty0 (2):\penalty0 980--994, 1994.

\bibitem[Papantoni-Kazakos and Gray(1979)]{PapantoniGray1979}
P.~Papantoni-Kazakos and R.~M. Gray.
\newblock Robustness of estimators on stationary observations.
\newblock \emph{The Annals of Probability}, 7\penalty0 (6):\penalty0 989--1002,
  1979.

\bibitem[Parthasarathy(1967)]{Parthasarathy1967}
K.~R. Parthasarathy.
\newblock \emph{Probability measures on metric spaces}, volume 352.
\newblock American Mathematical Soc., 1967.

\bibitem[Peligrad(1998)]{Peligrad1998}
M.~Peligrad.
\newblock On the blockwise bootstrap for empirical processes for stationary
  sequences.
\newblock \emph{Ann. Probab.}, 26\penalty0 (2):\penalty0 877--901, 1998.

\bibitem[Politis and Romano(1990)]{Politis1990}
D.~N. Politis and J.~P. Romano.
\newblock A circular block-resampling procedure for stationary data.
\newblock In \emph{Exploring the limits of bootstrap ({E}ast {L}ansing, {MI},
  1990)}, Wiley Ser. Probab. Math. Statist. Probab. Math. Statist., pages
  263--270. 1990.

\bibitem[Radulovi{\'c}(1996)]{Radulovic1996}
D.~Radulovi{\'c}.
\newblock The bootstrap for empirical processes based on stationary
  observations.
\newblock \emph{Stochastic Process. Appl.}, 65, 1996.

\bibitem[Rosenblatt(1956)]{Rosenblatt1956}
M.~Rosenblatt.
\newblock A central limit theorem and a strong mixing condition.
\newblock \emph{Proc. Nat. Acad. Sci. U. S. A.}, 42:\penalty0 43--47, 1956.

\bibitem[Sch\"{o}lkopf and Smola(2002)]{SchoelkopfSmola2002}
B.~Sch\"{o}lkopf and A.~J. Smola.
\newblock \emph{Learning with Kernels}.
\newblock Massachusetts Institute of Technology, Cambridge, 2002.

\bibitem[Shao and Yu(1993)]{ShaoYu1993}
Q.~M. Shao and H.~Yu.
\newblock Bootstrapping the sample means for stationary mixing sequences.
\newblock \emph{Stochastic Process. Appl.}, 48\penalty0 (1):\penalty0 175--190,
  1993.

\bibitem[Singh(1981)]{Kuensch1981}
K.~Singh.
\newblock On the asymptotic accuracy of {E}fron's bootstrap.
\newblock \emph{Ann. Statist.}, 9\penalty0 (6):\penalty0 1187--1195, 1981.

\bibitem[Steinwart and Christmann(2008)]{SteinwartChristmann2008}
I.~Steinwart and A.~Christmann.
\newblock \emph{Support vector machines}.
\newblock Information Science and Statistics. Springer, New York, 2008.

\bibitem[Steinwart et~al.(2009)Steinwart, Hush, and Scovel]{SteinwartHush2009}
I.~Steinwart, D.~Hush, and C.~Scovel.
\newblock Learning from dependent observations.
\newblock \emph{Journal of Multivariate Analysis}, 100:\penalty0 175--194,
  2009.

\bibitem[Straf(1969)]{Straf1969}
M.~L. Straf.
\newblock A general skorohod space, 1969.

\bibitem[Straf(1972)]{Straf1972}
M.~L. Straf.
\newblock Weak convergence of stochastic processes with several parameters.
\newblock pages 187--221, 1972.

\bibitem[Strohriegl and Hable(2016)]{StrohrieglHable2016}
K.~Strohriegl and R.~Hable.
\newblock On qualitative robustness for stochastic processes.
\newblock \emph{Metrika}, pages 895--917, 2016.

\bibitem[van~der Vaart(1998)]{VanderVaart1988}
A.~W. van~der Vaart.
\newblock \emph{Asymptotic statistics}, volume~3 of \emph{Cambridge Series in
  Statistical and Probabilistic Mathematics}.
\newblock Cambridge University Press, Cambridge, 1998.

\bibitem[Z\"ahle(2015)]{Zaehle2012}
H.~Z\"ahle.
\newblock Qualitative robustness of statistical functionals under strong
  mixing.
\newblock \emph{Bernoulli}, 21\penalty0 (3):\penalty0 1412--1434, 2015.

\bibitem[Z{\"a}hle(2016)]{Zaehle2016}
H.~Z{\"a}hle.
\newblock A definition of qualitative robustness for general point estimators,
  and examples.
\newblock \emph{Journal of Multivariate Analysis}, 143:\penalty0 12--31, 2016.

\end{thebibliography}

\end{document}